\documentclass[tbtags,reqno]{amsart}

\usepackage[latin1]{inputenc}
\usepackage{graphicx}
\usepackage{latexsym,amsmath,amssymb}
\usepackage{psfrag}
\usepackage{caption}

\newcommand{\Zm}{\mathbb{Z}}

\newcommand{\St}{\mathfrak S}

\newcommand{\lb}{[\![}
\newcommand{\rb}{]\!]}
\newcommand{\s}{\sigma}
\newcommand{\la}{\lambda}

\newcommand{\de}{\delta}
\newcommand{\eps}{\varepsilon}

\newcommand{\HI}{\mathcal{H}\mathcal{I}}
\newcommand{\HP}{\mathcal{H}\mathcal{P}}
\newcommand{\U}{\mathcal{U}}
\newcommand{\D}{\mathcal{D}}
\newcommand{\UU}{\mathbf{U}}
\newcommand{\DD}{\mathbf{D}}
\def\qee{$\hfill{\Box}$}

\newcommand{\K}{\mathbb{K}}

\newcommand{\Y}{\mathbf{Y}}

\newtheorem{theorem}{Theorem}
\newtheorem{proposition}[theorem]{Proposition}
\newtheorem{corollary}[theorem]{Corollary}
\newtheorem{lemma}[theorem]{Lemma}
\newtheorem{definition}[theorem]{Definition}

\author{Dominique Gouyou-Beauchamps}
\address{Universit\'e Paris Sud, 
Laboratoire de Recherche en Informatique, 91405 Orsay 
France}
\author{Philippe Nadeau}
\address{Fakult\"at f\"ur Mathematik, Universit\"at Wien, Nordbergstrasse 15, 1090 Vienna, Austria}

\title[Signed enumeration of ribbon tableaux]{Signed enumeration of ribbon tableaux}


\begin{document}

\begin{abstract}
We give an extension of the classical Schensted correspondence to the case of ribbon tableaux, where ribbons are allowed to be of different sizes. This is done by extending Fomin's growth diagram approach of the classical correspondence between permutations and pairs of standard tableaux of the same shape, in particular by allowing signs in the enumeration. As an application  we give a combinatorial proof for the column sums of the character table of the symmetric group. 
\end{abstract}

\maketitle

\section{Introduction}

The Schensted correspondence~\cite{Schen} is a bijection between permutations
and pairs of standard Young tableaux of the same shape. 
It has been extended in numerous
ways, the most famous being certainly the Robinson-Schensted-Knuth
correspondence between matrices of integers and pairs of semi-standard tableaux
of the same shape. Other extensions exist, for instance oscillating tableaux
\cite{Sundaram,ChauveDul,DulucqSag1,DulucqSag2}, skew tableaux~\cite{SaganStanley}, 
shifted Young tableaux~\cite{Sagan}, and $k$-ribbon tableaux~\cite{ShiWhi,StanWhi}.

Sergey Fomin developed a general theory of such correspondences, cf. 
\cite{FominGen,FominDual,FominSchen,FominSchur,FominStan}. It unifies
the correspondences listed above by interpreting these tableaux as paths in
so-called \emph{graded graphs}. For instance, a Young tableau in this context is
viewed as a particular kind of path in the graph whose vertices are integer
partitions, and where $(\la,\mu)$ is an edge if $\mu$ is a partition obtained
after adding $1$ to a part of $\la$. This graph is usually called the Young
graph (or Young lattice). Other combinatorial objects can be then represented by
considering other paths, for instance by modifying the extreme points of the
path, or the edges that one can use. This is a way of looking at oscillating or
skew tableaux inside the Young graph for instance. Then the local properties of
the graph will give rise to various bijective correspondences, all consequences
of one elementary bijection.

Furthermore, Fomin gives in parallel a linear algebraic approach to his
results, which is directly inspired by the work of Stanley
\cite{StanleyDiff}. As a matter of fact, most of Fomin's results
 have both a bijective and an algebraic proof.\medskip

In~\cite{WhiteBij}, Dennis White describes another extension of the
Schensted correspondence for {\em ribbon tableaux} where ribbons are allowed to
have 
different sizes; his goal was in fact to give a combinatorial proof of the
second orthogonality relation for characters of
 the symmetric group. The algorithm describing his correspondence is a
complicated insertion mechanism, along the lines of the original Schensted
correspondence; this forces him
moreover to put some restrictions on the ribbon tableaux he considers.\medskip

We will show here how the approaches of Fomin --both bijective
and algebraic-- can be adapted to apply to the correspondence of White.
There will be two benefits: first, this will extend the original work of White, 
by getting rid of his restrictions in particular.
Second, in the process we will have to extend Fomin's setting to graphs 
which are more general than the ones considered in~\cite{FominDual,FominSchen}.

Note that we will deal here with \emph{signed enumerations}, meaning that we will count objects with weights plus or minus one. It will appear that, for the bijective approach, we will have to appeal to the famous \emph{involution principle} of Garsia and Milne~\cite{GarMil}.
\medskip

The paper is organized as follows: in Section~\ref{sect:rib} we give the main definitions about ribbons and ribbon tableaux, as well as some elementary operations on these objects. In Section~\ref{sect:sign} we define different notions about signed enumeration, and we recall the famous Involution Principle of Garsia and Milne. Section~\ref{sect:hook} introduces hook permutations and involutions, which are the objects in correspondence with ribbon tableaux in our main results.  
Section~\ref{sect:results} states these results, which are Theorems
\ref{th:bij_perm} and~\ref{th:bij_invol}: there exists a signed bijection
between hook permutations and pairs of ribbon tableaux of the same shape.
 The description of the bijections is given in the two following sections: in Section~\ref{sect:local} we recall some local rules that depend on the operations introduced in Section \ref{sect:rib}, and Section~\ref{sect:bij} shows how to define a global correspondence from these local rules. The technical parts of the proofs are given in Appendix \ref{sect:appa}. In Section~\ref{sect:alg} an algebraic version of Theorem~\ref{th:bij_perm} is given. 
  
  Section~\ref{sect:char} contains an application of Theorem~\ref{th:bij_invol} to the column sums of the character table of the symmetric group. Finally, Section~\ref{sect:ext} explains how the methods developed can be used for other enumerations, and details in what ways this bijective and algebraic setting is a generalization of Sergey's Fomin graded graphs in duality. 

\section{Ribbons} \label{sect:rib}

\subsection{Definitions}

A \textit{partition} $\la=(\la_1,\ldots,\la_m)$ is a nonincreasing finite sequence of positive integers; these integers are the \textit{parts} of the partition, the \textit{size}
of the partitions being their sum $|\la|:=\sum_i \la_i$. A \textit{composition} $c$ is
as finite sequence of positive integers; we can associate to $c$
 a partition $\widetilde{c}$ by rearranging the integers in nonincreasing order.
A partition $\la$ of size $n$ can be described with the exponential notation
$\la=[1^{j_1},2^{j_2}\cdots n^{j_n}]$, where $j_i$ is the number of parts of size $i$.
If $j_i=0$, then $i^{j_i}$ is not written and if $j_i=1$, $j_i$ is not written.

We will identify a partition $(\la_1,\ldots,\la_m)$ with its Ferrers diagram,
which is the left justified set 
of cells (i.e. unit squares of $\Zm^2$) such that the $i$th row from the top
contains $\la_i$ cells; the diagram on the left of Figure $1$ represents the
partition $(8,6,5,2,2,1,1)=[1^2\,2^2\,5\,6\,8]$.

\begin{figure}[ht]
\centering
\includegraphics[width=0.8\textwidth]{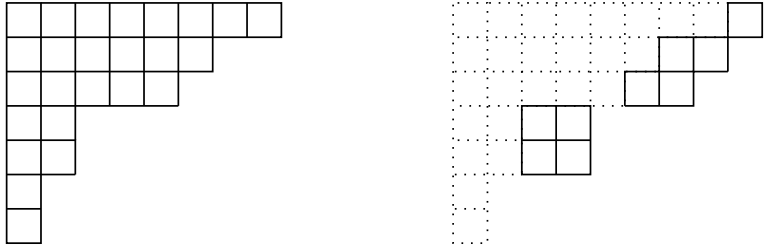}
\caption{Partition and skew shape.\label{fig:Shapes}}
\end{figure}

Let $\Y$ be the set of integer partitions, and $\Y_n$ the subset of partitions
of size  $n$. Two partitions $\la\subseteq \mu$ (in the sense of inclusion of
Ferrers diagrams) define a  {\em skew shape} $\mu/\la$. We will identify here a
skew shape with the set of cells $\mu\backslash \la$,  whenever $\mu$ or $\la$
is clear from the context; though, in general, two distinct skew shapes may
define the same set, this will not create any ambiguity here. The {\em size} of
$\mu/\la$ is its number of cells $|\mu|-|\la|$ and will be noted $|\mu/\la|$.
The skew shape $(9,8,7,4,4,1,1)/(8,6,5,2,2,1,1)$  represented on the right of
Figure~\ref{fig:Shapes} has size $9$.

 Let us say that a subset $S$ of cells of  $\Zm^2$ is {\em connected} if, for
every two cells $c,c'$ in $S$, there exists a sequence of cells
$c=c_0,c_1,\ldots,c_t=c'$ in $S$ such that $c_i,c_{i+1}$ have a common side for
all $i$. We can then define the notion of ribbons:

\begin{definition}
A \emph{ribbon} is a connected skew shape that does not contain a $2$ by $2$
square of cells.
\end{definition}

\begin{figure}[t]
\centering
\includegraphics[width=0.7\textwidth]{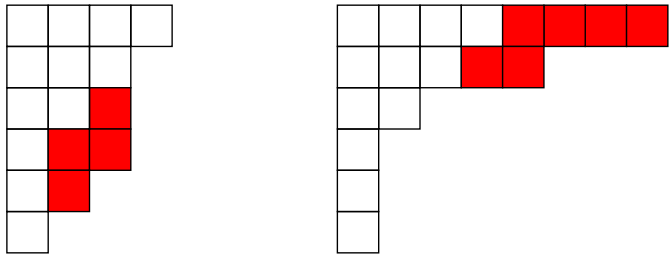}
\caption{Examples of ribbons.\label{fig:Ribbons}}
\end{figure}

Let $r=\mu/\la$ be a non empty ribbon. Such a ribbon is said to be
 $\mu$-addable and $\la$-removable. Its {\em height}
 $h(r)$ is defined as the number of lines it occupies  minus
 one, and its {\em sign} is then $\eps(r):=(-1)^{h(r)}$. By convention
we will set $\eps(\la/\la):=1$. The bottom left cell of $r$ is its {\em tail},
 and the top right one is its {\em head}. Given a partition $\la$, the ribbons
 that can be removed or added to $\la$ are entirely determined by the
coordinates of their heads and tails. On Figure~\ref{fig:Ribbons} are
two ribbons: the left one has size $4$, height $2$ and sign $+1$, and
the right one has size $6$, height $1$ and thus sign $-1$.

A {\em hook} is a non empty ribbon of shape ${\lambda/\emptyset}$,
which is equivalently a partition of the kind $(k,1,\ldots,1)$. Note
that a hook is characterized by the data of its size $s$ and
height $h\in \lb 0,s-1\rb$.

\medskip

\begin{figure}[!ht]
\centering
\includegraphics[width=0.8\textwidth]{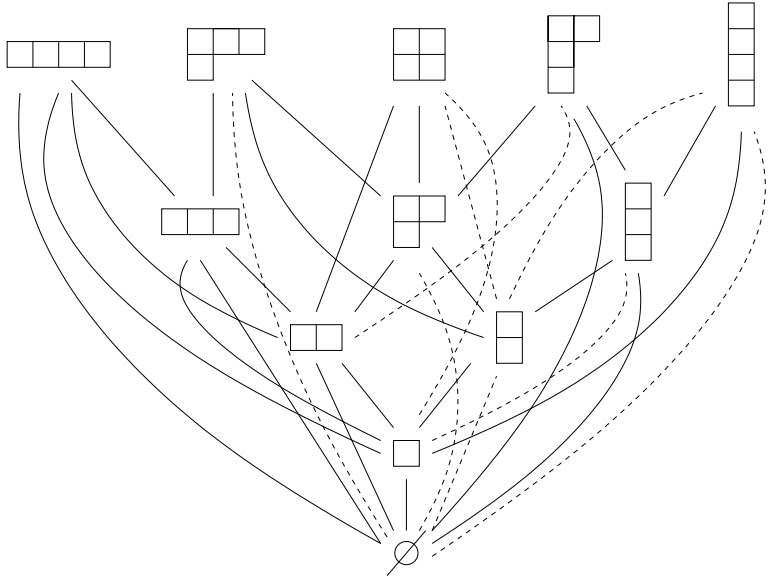}
\caption{First levels of the ribbon graph $GR$.\label{fig:Ribbon_Graph}}
\end{figure}

For $i$ a positive integer, we note $Rib_i$ the set of ribbons of
size $i$, and $Rib:=\bigcup_{i>0} Rib_i$. One should think of $\Y$ as
 the vertices of a graph $GR$ whose edges are the elements of $Rib$;
 each edge carries in addition a sign which is simply the sign of
 the corresponding ribbon. $GR$ is structured in different levels
 given by the partitions of a given size. Figure~\ref{fig:Ribbon_Graph}
 shows the first levels of $GR$, where the dotted edges correspond
to negative ribbons. Adding (respectively removing) a ribbon corresponds
 simply to making a step up (resp. down) in $GR$.

\medskip

\begin{definition}
 A {\em ribbon tableau} of shape $\la\in \Y$
and length $\ell$ is a sequence of partitions  $\la^0=\emptyset\subset
\la^1\subset
\ldots\la^{\ell-1}\subset \la^\ell=\la$ such that $r_i:=\la^{i+1}/\la^i$
 is a nonempty ribbon for every $i<\ell$.
\end{definition}
\medskip

We will often represent a ribbon tableau by a labeling of the cells of
 $\la$, in which the cells labeled $i$ coincide with the ribbon $r_i$.
 Note that a ribbon tableau is equivalently a path of length $\ell$
 in the graph $GR$, going up from  $\emptyset$ to $\la$; this
 interpretation is the key to the extensions described in Section
\ref{sect:ext}.

We need some more definitions. 
The \emph{sign} $\eps(P)$ of a tableau
$P$ is the product of the signs of the ribbons $r^{(i)}$.
The \emph{content} $c(P)$  is the composition of $|\la|$
in $\ell$ parts formed by the sequence of sizes $|r^{(i)}|$.
We will note $RT_{\la,c}$ the tableaux of shape $\la$ and
content $c$, where $c$ is a given composition of $|\la|$ and
$RT_{\la,l}$ the set of all the ribbon tableaux of shape $\la$ and
length $l$.
\medskip

\begin{figure}[ht]
\centering
\includegraphics[width=0.5\textwidth]{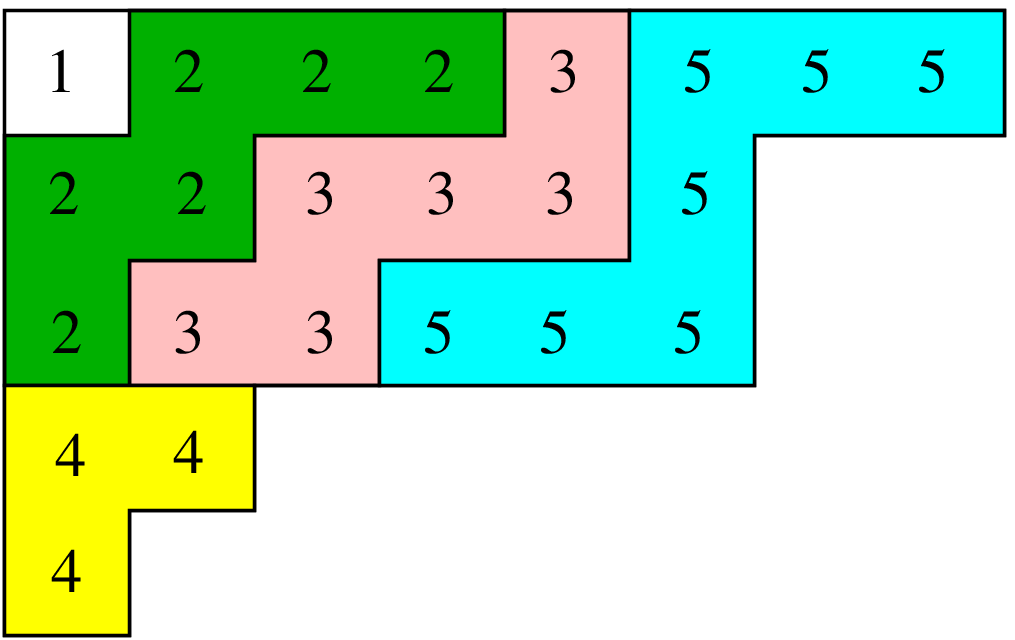}
\includegraphics[width=0.8\textwidth]{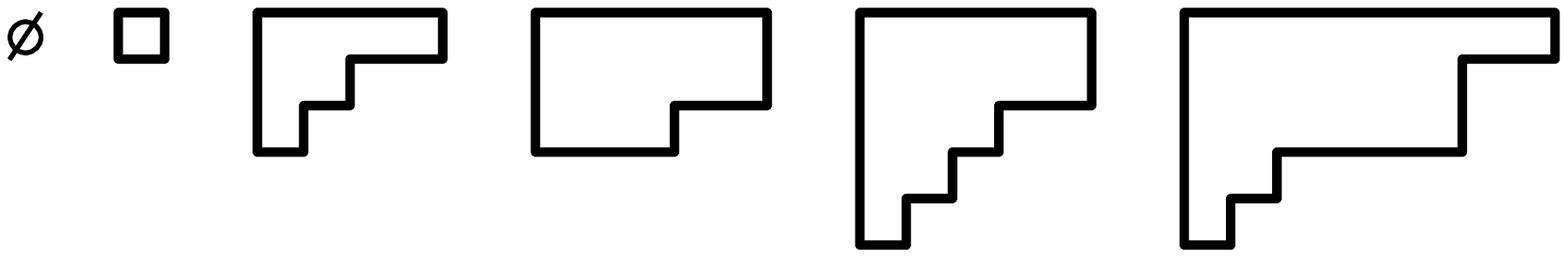}
\caption{\label{fig:Tableau} A ribbon tableau: two equivalent representations.}
\end{figure}

Figure~\ref{fig:Tableau} shows a tableau of
shape $(8,6,6,2,1)$, content $(1,6,6,3,7)$ and sign
$(-1)^0(-1)^2(-1)^2(-1)^1(-1)^2=-1$.

\medskip

\subsection{Operations on ribbons} \label{subsect:oper} We will now introduce
 some classical operations on ribbons which are necessary for the
definition of the {\em local rules} of Section~\ref{sect:local}.
\medskip

\begin{figure}[!ht]
\includegraphics[width=\textwidth]{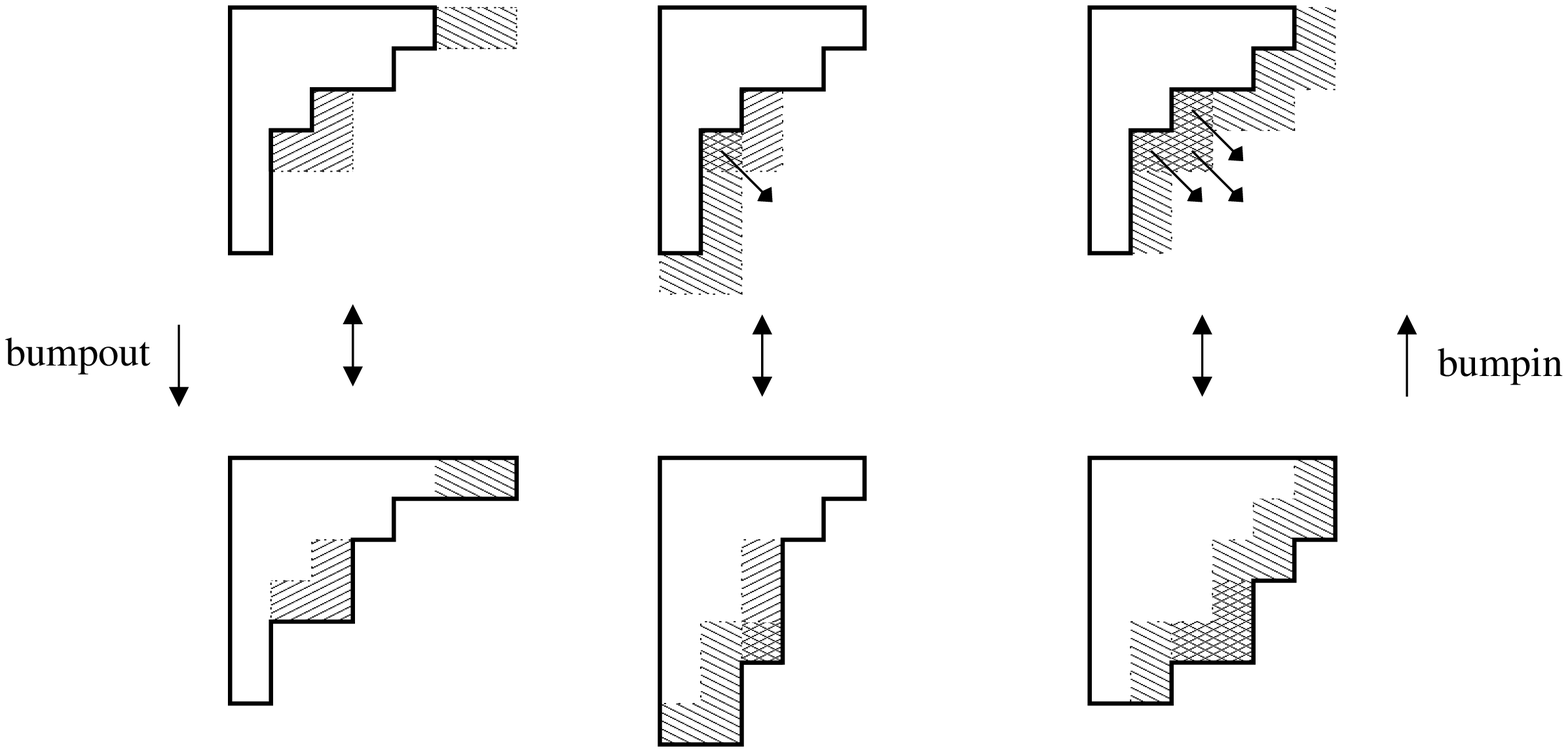}
\caption{\label{fig:bump} Operations $bumpin$ and
$bumpout$ }
\end{figure}

$\mathbf{bumpin,~bumpout}$: let $\la$ be a partition, and $r,r'$ be two ribbons
that are $\la$-addable, such that $r$ and $r'$ have different
 heads and different tails. Then $bumpout(r,r')$ is the set of cells
 $(r\backslash r')\cup (r'\backslash r)\cup (r'\cap
r)_{\searrow}$ where $A_{\searrow}$ is the translated of $A$ by
the vector $(1,-1)$. The operation $bumpin(r,r')$ is similarly defined  for two $\xi$-removable ribbons, by translating the common cells between $r$ and $r'$ by the vector
$(-1,1)$ (these definitions vary slightly from the ones commonly used). Figure~\ref{fig:bump} shows the way this operation acts, according to the relative positions of $r$ and $r'$: they can be disjoint, or partially overlap, or one can be included in the other.

$\mathbf{prev,next,first}:$ Let $\la$ be a partition, $k$ a positive integer, and $h$ a nonnegative integer. A result of Shimozono and White~\cite{ShiWhi} is that if $(r_i)_{i=0\ldots t}$ (respectively $(r'_i)_{i=1\ldots t'}$) are all ribbons of size $k$ and height $h$ that are $\la$-addable (resp. $\la$-removable), then (1) t=t', and (2) the enumeration order of the ribbons can be chosen so that $r_0<r'_1<r_1<\ldots <r'_t<r_t$ where $rib_1<rib_2$ means that the head of $rib_1$ is north east of the head of $rib_2$.

\begin{figure}[!ht]
\centering
\includegraphics[width=0.3\textwidth]{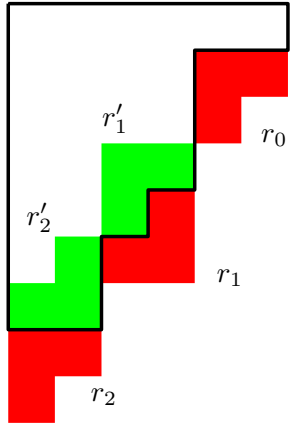}
\caption{\label{fig:first} Addable and removable ribbons of height $1$ and size $3$.} 
\end{figure}
 
Figure~\ref{fig:first} illustrates Shimozono and  White's result. This allows to define certain operations $first,next$ and $prev$ on ribbons:

\begin{itemize}

\item if $hk$ is a hook of size $k$ and height $h$, we define 
$first(\la,hk)$ as the ribbon $r_0$ above.

\item If $r'$ is a $\la$-removable ribbon, i.e. $r'=r'_i$ for a
 certain  $i\in\lb 1,t\rb$, then $next(\la,r'):=r_i$.

\item Reciprocally, if $r=r_i$ for any $i\in \lb 1,t\rb$ 
is a $\lambda$-addable ribbon, we define $prev(\la,r_i):=r'_i$ for
$i\geq 1$, and $prev(\la,r_0):=\emptyset$.
\end{itemize}

\begin{figure}[!ht]
\includegraphics[width=\textwidth]{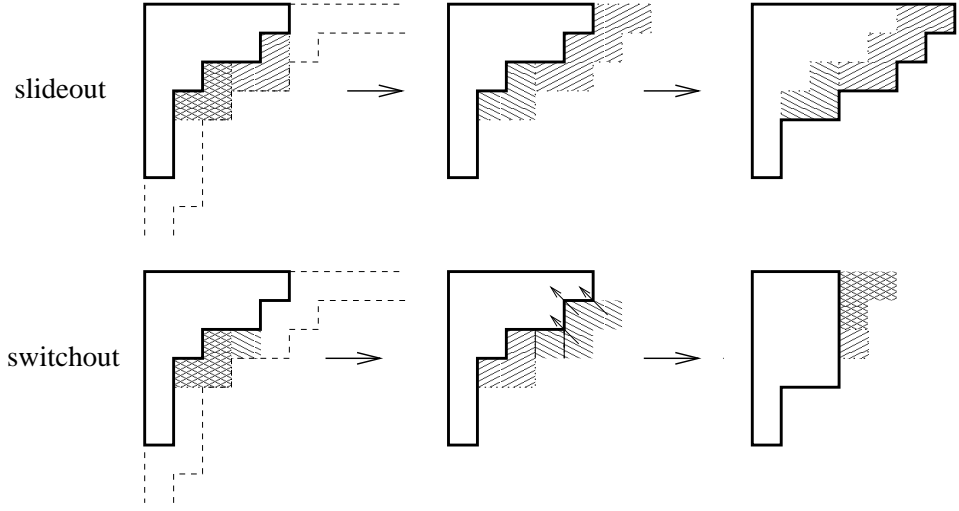}
\psfrag{swichout}{\Large $switchout$} \psfrag{slideout}{\large
$slideout$} \caption{Operations $switchout$ and $slideout$.
\label{fig:switch}}
\end{figure}

$\mathbf{slideout,switchout,slidein,switchin}$: let $\la$ be a partition, and $r,r'$ two $\la$-addable ribbons
, with identical tails but different heads; we assume without loss of generality that $|r|>|r'|$. The {\em external band} of $\la$ consists 
of all cells between the infinite south east boundary of $\la$ and its translated by $(1,-1)$, i.e. the cells enclosed by the dotted line on Figure~\ref{fig:switch}. Let $\tau$ be the subset of the external band formed by the $|r'|$ connected cells, north west of $r$ and adjacent to it. Then two cases can occur:
\begin{itemize}
\item if $\tau \cup r$ forms a $\la$-addable ribbon, we define the partition $slideout(\la,r,r')=\la\cup r\cup \tau$.

\item otherwise, we define $switchout(\la,r,r')=(\la\cup r')\backslash \tau_{\nwarrow}$, where  $\tau_{\nwarrow}$ is
 the translated of $\tau$ by the vector $(-1,1)$.

\end{itemize}

If $r$ and $r'$ have the same head but different tails,
one performs the same operations on the transposed partitions.
The operations $switchin$ and $slidein$ are defined similarly
 on $\la$-removable ribbons: see White~\cite{WhiteBij} for supplementary explanations.

\section{Signed sets and signed bijections} \label{sect:sign}

In this work we have to deal with signed enumerations, so we need
some definitions and notations to explain what we mean by a
bijection in this context. All sets are assumed to be finite.

\begin{definition}[Signed Sets]
A signed set is a set $A$ together with a decomposition $A=A^+\cup
A^-$ where $A^+\cap A^-=\emptyset$. 
The members of $A^+$ are \emph{positive} elements, those of
$A^-$ are \emph{negative}.
\end{definition}

Such a decomposition is equivalent to a function  $\delta:A\rightarrow
\{1,-1\}$, with the obvious
correspondence $A^+=\delta^{-1}(\{1\})$ and
$A^-=\delta^{-1}(\{-1\})$. Our objects of study here are the sets
 $RT_{\lambda,\mu}$ , the sign being given by the function $\eps$.
  Note also that unless
explicitly stated, usual sets are considered as positive sets.

A function $f$ between two signed sets is \emph{sign preserving} (resp.
\emph{sign reversing}) if $a$ and $f(a)$ have the same sign (resp.
different signs) for all $a$. Fixed points of a function $i$ form
the set $Fix(i)$.

\begin{definition}[Signed bijections]
A signed bijection between the signed sets $A$ and $B$ is the data
of 3 functions $i_A,i_B$ and $\phi$ such that $i_A$ (resp. $i_B$)
is an involution on $A$ (resp. $B$) which is sign reversing
 outside of its fixed point set, and $\phi$ is a
sign preserving bijection between $Fix(i_A)$ and $Fix(i_B)$.
\end{definition}

The {\em signed cardinal} (or {\em signed sum}) of a signed set $A$ is $|A|_\pm
=|A^+|-|A^-|$; if the sign is given by a function $\delta$, then
 we have $|A|_\pm=\sum_{a\in A} \delta(a)$. A signed bijection
  between $A$ and $B$ proves that $|A|_\pm=|B|_\pm$; indeed we have
$$|A|_\pm=|Fix(i_A)|_\pm=|Fix(i_B)|_\pm=|B|_\pm.$$

The central equality comes from the sign preserving bijection $\phi$, the other
ones from the fact that $i_A$ and $i_B$ are sign reversing, so the pairs
$\{a,i_A(a)\}$ such that $a\neq i_A(a)$ have a zero contribution to the signed
cardinal of $A$, the same being true for $B$ and $i_B$. Now $|A|_\pm=|B|_\pm$ is
equivalent to $|A^+|+|B^-|=|B^+|+|A^-|$, and a bijection proving this equality
(i.e. a bijection $\psi$ between the usual sets $A^+\sqcup B^-$ and
$B^+\sqcup A^-$) is clearly equivalent to a signed bijection between $A$ and
$B$. This explains why signed bijections are the correct generalizations of
bijections, in that they give a combinatorial explanation of the equality of
signed cardinals. \medskip

\textbf{The Involution Principle of Garsia and Milne}\smallskip

Garsia and Milne gave the first fully bijective proof of a combinatorial version
of a famous identity of Rogers-Ramanujan~\cite{GarMilshort,GarMil}: this states
that the number of partitions $(\la_1,\ldots,\la_k)$ of $n$ verifying
$\la_i-\la_{i+1}\geq 2$ for all $i<k$ is the same as those verifying
$\la_i\equiv 1\text{ or }4$ modulo $5$ for all $i\leq k$. To achieve this, they
defined and used a general principle, that we recall now.\medskip

Let $A,B$ be two finite signed sets. Let also $i_A, i_B$
 be two involutions on $A$ and $B$ respectively, and
$\phi$ a bijection between $A$ and $B$. We
suppose that $\phi$ preserves signs, whereas $i_A$ and $i_B$
reverse signs outside their fixed point sets.

Under those assumptions one has clearly $|Fix(i_A)|_\pm=|Fix(i_B)|_\pm$, but not
necessarily a {\em signed bijection} proving this equality. The principle of
Garsia and Milne is the construction of such a signed bijection $(\psi,j_A,j_B)$
between $Fix(i_A)$ and $Fix(i_B)$ in the following manner: let $a\in A$. We
apply to it the function $\phi:A\rightarrow B$, then alternatively, the
functions $\phi^{-1}\circ i_B:B\rightarrow A$ and $\phi\circ i_A:A\rightarrow
B$, stopping as soon as the element $x$ obtained is: \begin{itemize} \item
either in $Fix(i_A)$, in which case one sets $j_A(a):=x$; \item or in
$Fix(i_B)$, in which case one sets $\psi(a):=x$ (and $j_A(a):=a$). \end{itemize}

To define $j_B$ (and $\psi^{-1}$), one uses the symmetric
procedure starting from $b\in B$. These procedures terminate and
 give the seeked signed bijection: see for instance~\cite[p.76]{Kerber} for a
proof.

\section{Hook permutations and hook involutions} \label{sect:hook}

We now introduce the notions of \textit{hook permutations and hook involutions},
which play the role of the ordinary permutations and involutions of the
Schensted correspondence in our main results stated in Section~\ref{sect:results}.

\subsection{Hook Permutations}

If $H=(H_1,\ldots,H_\ell)$ is an ordered sequence of $\ell$ hooks,
its {\em content} $c(H)$ is the composition $(|H_1|,\ldots,|H_\ell|)$.

\begin{definition}
A hook permutation $(H,\s)$ is an ordered sequence
$H=(H_1,\ldots,H_\ell)$ of $\ell$ hooks, together with a
permutation $\s$ of $\lb \ell\rb$. The length of a hook
permutation is $\ell$, its size is $\sum_i |H_i|$, and its content
$c(H,\s)$ is the composition $(|H_{\s(1)}|,\ldots,|H_{\s(\ell)}|)$.
\end{definition}

\begin{figure}[ht]
\centering
\includegraphics[width=\textwidth]{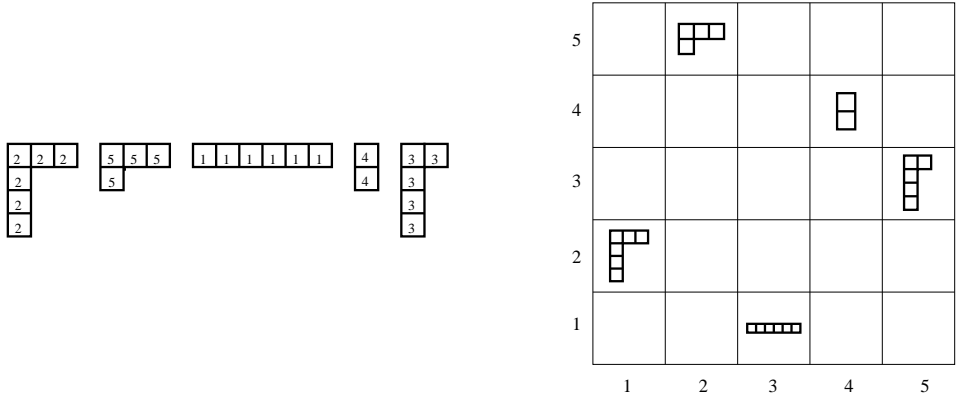}
\caption{Two representations of the same hook permutation of length
$5$, size $23$ and content $(6,4,6,2,5)$. \label{fig:Hook_Perm}}
\end{figure}

We will write $\HP$ for the set of hook permutations, its elements
of content $\mu$ forming $\HP(\mu)$ (where $\mu$ is any
composition). Hook permutations can be represented by the list $H$
where the cells of hook $H_i$ are numbered by $\s(i)$, or by
square arrays of size $\ell$ such that entry $(i,j)$ is empty
unless $j=\s(i)$ in which case it is occupied by the $i$th hook
$H_i$. Illustrations are given Figure~\ref{fig:Hook_Perm}.

We then have the following proposition:

\begin{proposition}
\label{prop:enumhookperm}
The number of hook sequences of length $\ell$ and total size $n$
is equal to $\binom{n+\ell-1}{2\ell-1}$.
\end{proposition}

\begin{figure}[!ht]
\centering
\includegraphics[width=0.7\textwidth]{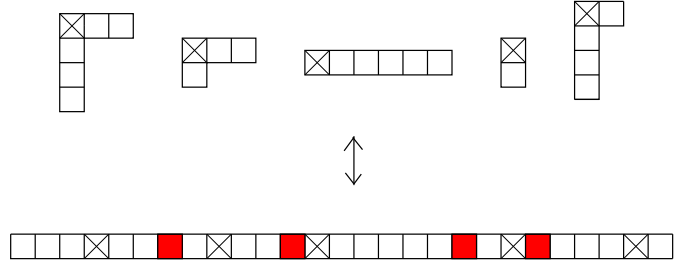}
\caption{The bijection in Proposition~\ref{prop:enumhookperm}.
\label{fig:hook_list}}
\end{figure}

\noindent{\bf{Proof:~}} Given such a hook sequence, we can associate 
  to it $2\ell-1$ integers in $\lb n+\ell-1\rb$, as illustrated
graphically on Figure~\ref{fig:hook_list}, in which $n=23$,
$\ell=5$, so $n+\ell-1=27$ and $2\ell-1=9$, and the subset
 of integers is $\{4,7,9,12,13,19,21,22,26\}$.

Conversely, if we have $2\ell-1$ integers
$1\leq i_1<\ldots<i_{2\ell-1}\leq n+\ell-1$, and we set by convention
 $i_0=0$ and $i_{2\ell}=n+\ell$, then the list of hooks
$(H_1,\ldots,H_\ell)$ in bijection is characterized by the
fact that $H_i$ is the hook of size $i_{2i}-i_{2i-2}-1$ and height
$i_{2i-1}-i_{2i-2}-1$. This is clearly bijective.
\qee

\subsection{Hook Involutions}

\begin{definition}
 Hook involutions are hook permutations such that their array
representation is symmetric with respect to the diagonal $i=j$.
\end{definition}

In other words, these are hook permutations $(H,\s)$ such that $\s$
is an involution, and $H_i=H_j$ if $j=\s(i)$. For a hook
involution $I=(H,\s)$, we define its {\em sign} as
 $\eps(I)=\prod_{i/\s(i)=i} \eps(H_i)$. It is the product of the
signs of the hooks associated to fixed points.

 We note $\HI$ the signed set of hook involutions,
 $\HI(\mu)$ the signed subset of those hook involutions with content
 $\mu$, and finally ${\HI_{spec}(\mu)}$ the elements of
 $\HI(\mu)$ all of whose fixed points are hooks of odd size and of
 height $0$. Note that all elements of ${\HI_{spec}(\mu)}$ are then positive.

\begin{lemma} \label{lem:inv_hook_inv}
There is a sign reversing involution on
 $\HI(\mu)\backslash {\HI_{spec}(\mu)}$.
\end{lemma}

\noindent {\bf Proof: } If $I=(H,\s)\notin {\HI_{spec}(\mu)}$
then let $i=\s(i)$ be its smallest fixed point contradicting
the definition of $\HI_{spec}(\mu)$. Let $h$ be the height of
 the hook $H_i$. If $H_i$ is of even size , then we let $H'_i$
 be the hook of the same size and of height $h+1$ (resp. $h-1$)
 if $h$ is even (resp. odd). If $H_i$ is of odd size, so that
 necessarily $h\neq 0$ by the definition of $H_i$, then we let
 $H'_i$ be the hook of the same size and of height $h+1$
 (resp. $h-1$) if $h$ is odd (resp.even).

 Let $H'$ be the hook list equal to $H$ except in position $i$
where $H'_i$ replaces $H_i$. If we define $f(I)=(H',\s)$, then we
have the desired sign reversing involution on $\HI(\mu)\setminus
\HI_{spec}(\mu)$. \qee

\begin{corollary} \label{corol:inv_hook_inv}
For any composition $\mu$, $|\HI(\mu)|_\pm=|{\HI_{spec}(\mu)}|.$
\end{corollary}

We will give some consequences of this result in Section~\ref{sect:char}.

\section{Main results} \label{sect:results}

The first result is a generalization of the
Schensted correspondence. We note by $Id$ the identity function on hook
permutations.

\begin{theorem} \label{th:bij_perm}
  Let $n,\ell$ be two positive integers. There exists an explicit signed
bijection $(Id,i,\phi)$ between hook permutations of size $n$ and length $\ell$,
and pairs of ribbon tableaux of size $n$ and length $\ell$ with the same shape.

   This bijection {\em  preserves contents}, which means: if
$i(P,Q)=(P_1,Q_1)$, then $c(P)=c(P_1)$ and $c(Q)=c(Q_1)$; and if
$\phi(\s,H)=(P,Q)$, then $c(H)=c(Q)$ and  
$c(H,\sigma )=c(P)$.
\end{theorem}

We will define the signed bijection is given in Section~\ref{sect:bij}, based on local rules given in Section~\ref{sect:local}; the proof of the correctness of the bijection is in Appendix~\ref{sect:appa}. In the special case where contents are partitions with certain constraints, then the preceding theorem is equivalent to the main result of White in~\cite{WhiteBij}. The idea
here is to use Fomin's techniques~\cite{FominDual,FominSchen} in the proof of this result: this sheds a new light on White's result, and lends itself to generalization in a more straightforward fashion.

The theorem has the following consequences concerning the
 signed enumeration of ribbon tableaux:

\begin{corollary} \label{corol:enum_pairs}
Let $\mu,\nu$ be two compositions of $n$ with $\ell$ parts, and write $\widetilde{\mu}=[1^{j_1},2^{j_2}\cdots]$. Then
\begin{align}
 \sum_{\substack{\la\in \Y_n\\ P\in RT_{\la,\mu},Q
\in RT_{\la,\nu}}} \varepsilon(P)\varepsilon(Q) &=
\delta_{\widetilde{\mu}\widetilde{\nu}}\cdot
1^{j_1}(j_1!)2^{j_2}(j_2!)\ldots ;\\
\sum_{\substack{\la\in \Y_n\\ P,Q\in RT_{\la,\ell}}}
    \varepsilon(P)\varepsilon(Q) &=
   \binom{n+\ell-1}{2\ell-1}\cdot \ell!
\end{align}
\end{corollary}

We will show that this corollary can be proved by techniques of linear algebra in Section~\ref{sect:alg}. Finally, the Schensted correspondence has the property that it restricts to a bijection between involutions of $\St_n$ and standard tableaux of size $n$. We will prove a version of this result for ribbon tableaux:

\begin{theorem} \label{th:bij_invol}
Let $n,\ell$ be two positive integers. There exists an explicit signed bijection between hook involutions of size $n$ and length $\ell$, and ribbon tableaux of size $n$ and length $\ell$.\end{theorem}

As explained in Section~\ref{sect:bij}, this cannot be deduced from White's correspondence for pairs of tableaux. We will prove this Theorem in Section~\ref{sect:bij} and Appendix~\ref{sect:appa}.



\section{Local Rules} \label{sect:local}

We wish to extend local the rules used by Shimozono and White~\cite{ShiWhi} to deal with ribbons of all possible sizes; this will be done simply by reformulating White's insertion rules of~\cite{WhiteBij} as local rules.
\emph{In all this section $\mu,\nu$ are two partitions of
respective sizes $m$ and $n$.}

\medskip

 For $i$ a nonnegative integer, we define $\U_i(\mu,\nu)$ as the
 set of partitions of size $\max(m,n)+i$ such that $\xi/\mu$
 and $\xi/\nu$ are ribbons.
$\U_i(\mu,\nu)$ is a signed set through $sgn(\xi):=
\eps(\xi/\mu)\cdot \eps(\xi/\nu)$. Similarly, we
define $\D_i(\mu,\nu)$ the set of partitions of size $\min(m,n)-i$
 such that $\la/\mu$ and $\la/\nu$ are ribbons; $\D_i(\mu,\nu)$ is
 a signed set through $sgn(\la):= \eps(\mu/\la)\cdot \eps(\nu/\la)$.
 ( Note that, as signed sets, $\U_i(\mu,\mu)$ and
$\D_i(\mu,\mu)$ contain only positive elements.)\medskip

\noindent
   \begin{minipage}[b]{0.7\linewidth}
\hspace{0.3cm} We draw a square where $2$ corners are labeled by
$\mu$ and $\nu$ as shown on the right. The bottom left corner will be
 labeled by partitions $\la$ from a set $\D_i(\mu,\nu)$, the top right one
 by partitions $\xi$ from a set $\U_i(\mu,\nu)$.
 In the case $\la=\mu=\nu$, the interior $C$ may
 be marked by a {\em nonempty hook}, or be left empty;
in all other cases it is left empty.
         \end{minipage}\hfill
   \begin{minipage}[b]{0.3\linewidth}
      \psfrag{L}{$\la$}
      \psfrag{M}{$\mu$}
      \psfrag{N}{$\nu$}
      \psfrag{h}{$C$}
      \psfrag{X}{$\xi$}
      \centering \includegraphics[width=0.5\textwidth]{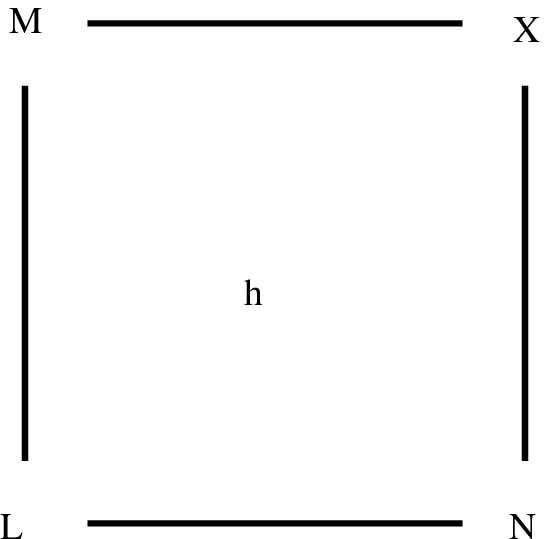}
   \end{minipage}

\medskip

To define local rules, it is necessary to use the
operations on ribbons and partitions defined in Section~\ref{sect:rib}.

Let $((\la,C),\mu,\nu)$ be given in a square, as above.
 What we mean by  {\em applying a local rule} is the following:
first find out which case applies in the list below, then {\em erase
$\la$ and $C$ from the square}, and finally write the outcome
 of the rule in the appropriate corner: the top right one for
 the rules D1 to D6, and the bottom left one for the rule S.

\smallskip

{\bf Direct rules: } in the following rules $\la$
is an element of a certain $D_i(\mu,\nu)$, $C$ is empty
 unless possibly when $\la=\mu=\nu$ in which case it may
 be filled by a hook. Let $r$, $r'$ be the ribbons $\mu/\la$
 and $\nu/\la$ (allowing empty ribbons).

\begin{itemize}

\item[$\bullet$] If  $\la=\mu=\nu$ and $C$ is empty, then
$\xi:=\la$.~~(D1)

 \item[$\bullet$] If  $\la=\mu=\nu$ and $C$ is a nonempty
 hook $eq$, then $\xi:=\la\cup first(\la,eq)$.~~(D2)

 \item[$\bullet$] If
$\la\neq\mu=\nu$, then $\xi:=\mu\cup next(\mu,\mu/\la)$.~~(D3)

\item[$\bullet$] If $\la=\mu\neq\nu$ (\emph{resp.}
$\la=\nu\neq\mu$), then $\xi:=\nu$ (\emph{resp.}
$\xi:=\mu$).~~(D4)

\item[$\bullet$] If  $\la\neq \mu \neq \nu$, then:
   \begin{itemize}
       \item if $r$ and $r'$ have neither the same tail nor the same head, 
       then $\xi:=\la\cup bumpout(r,r')$.~~(D5)
       \item if $r$ and $r'$ have the same head but different tails,
        or the same tail but different heads, then:
       \begin{itemize}
             \item if $slideout(\la,r,r')$ is well defined, then
             $\xi:=slideout(\la,r,r')$.~~(D6)
             \item otherwise, we set $\widehat{\la}:=
             switchout(\la,r,r')\in D_i(\mu,\nu)$.~~(S)
       \end{itemize}
   \end{itemize}
\end{itemize}

\smallskip

{\bf Inverse rules: } Now $\xi$ belongs to a certain set $U_i(\mu,\nu)$. We
write $r$,
$r'$ for the ribbons $\xi/\mu$ and $\xi/\nu$ (possibly empty).
$C$ is left empty except in rule (I2).

\begin{itemize}
\item[$\bullet$] If $\xi=\mu=\nu$ , then $\la:=\xi$.~~(I1)

\item[$\bullet$] If  $\xi\neq\mu=\nu$, then
\begin{itemize}
\item If $prev(\xi,r)=\emptyset$, we define $\la:=\mu$ and
$C$ is filled with the unique hook of size $|r|$ and height
$h(r)$;~~(I2)
 \item otherwise $\la:=\mu\backslash prev(\xi,r)$.~~(I3)
\end{itemize}

\item[$\bullet$]  If  $\xi=\mu\neq\nu$ (\emph{resp.}
 $\xi=\nu\neq\mu$), then $\la:=\nu$ (\emph{resp. }$
\la:=\mu$). ~~(I4)

\item[$\bullet$] If  $\xi\neq \mu \neq \nu$, then:
   \begin{itemize}
       \item if $r$ and $r'$ have neither the same tail nor the same head, then 
$\la:=\xi\backslash
       bumpin(r,r')$.~~(I5)
       \item if $r$ and $r'$ have the same head but different tails,
        or the same tail but different heads, then:
       \begin{itemize}
             \item if $slidein(\xi,r,r')$ is defined, then
             $\la:=slidein(\xi,r,r')$;~~(I6)
         \item otherwise we set $\widehat{\xi}:= switchin(\xi,r,r')\in
U_i(\mu,\nu)$.~~(T)
       \end{itemize}
   \end{itemize}
\end{itemize}

\medskip

\begin{proposition} \label{prop:locrul}
The rules D1 to D6 are the respective inverses of I1 to I6;
 S and T are involutions. Furthermore, D1-D6 and I1-I6
preserve signs between $\D_i(\mu,\nu)$ and $\U_i(\mu,\nu)$,
whereas S and T are sign reversing on
$\D_i(\mu,\nu)$ and $\U_i(\mu,\nu)$ respectively.
\end{proposition}

{\noindent{\bf Proof:~}} All these properties are already proved elsewhere,
albeit sometimes in a different form. For D2,D3 and I2,I3 this was proved by
Shimozono and White~\cite{ShiWhi}. For the rules D5,D6,I5,I6, S and T, the
result can be found in White~\cite{WhiteBij}. We will anyway give the proof for
the rule $S$ in Appendix~\ref{sect:appb}, using the encoding of partitions by infinite sequences $\de(\la)$:
we wish to show how this encoding is particularly suited to the study of ribbons.\qee
\medskip

Let us sum up the local rules in terms of signed bijections, since this will be useful in particular in the algebraic approach of Section~\ref{sect:alg}:

\begin{proposition}

\label{th:loc_bij}
Let $\mu,\nu$ be two partitions and $i$ a positive integer.
\begin{itemize} 
 \item[(a)] There exists an explicit bijection $\phi_1$ between $\U_i(\mu,\mu)$
and
  $\D_i(\mu,\mu)\sqcup \lb 0,i-1 \rb$. \item[(b)] For $\mu\neq
\nu$, there exists a {\em signed} bijection $(i_D,i_U,\phi_2)$ between 
$\D_i(\mu,\nu)$ and $\U_i(\mu,\nu)$.
\end{itemize}
\end{proposition}

{\noindent{\bf Proof:~}}  (a) is given by rules D2 and D3; for (b),
 the involutions $i_D$ and $i_U$ are given by rules S and T respectively,
 and the bijection $\phi_2$ consists of rules D4,D5 and D6.
\qee\medskip

This proposition has to be interpreted as a {\em local property }
of the graph $GR$: given $\mu$ and $\nu$ of size $m$ and $n$, it gives a
combinatorial link between vertices adjacent to both $\mu$ and $\nu$ at the
level $\min(m,n)-i$, and at the level $\max(m,n)-i$. Sections~\ref{sect:bij} and
\ref{sect:alg} will use this local structure to deduce global results on ribbon
tableaux. 

\section{Bijective approach} \label{sect:bij}

We want to use {\em growth diagram} techniques (see~\cite{FominSchen}) to carry out a
signed correspondence proving Theorems~\ref{th:bij_perm} and~\ref{th:bij_invol}.
This will be done in this section, but requires more work than a simple
application of Fomin's setting. In order to prove Theorem~\ref{th:bij_perm}, we
will actually need to make some {\em back and forth} moves in a growth diagram;
the correction of the correspondence will rely on the Involution Principle.

\subsection{Bijection for hook permutations}
\label{sub:bijhookperm}

Fix a positive integer $\ell$, and let $G_\ell$ be the grid of size $\ell\times
\ell$, made of $\ell^2$ squares (and $(l+1)^2$ vertices). Square $(i,j)$ is at
the intersection of the $i$th column from the left and the $j$th row from the
bottom of the grid. We wish now to label some of the vertices by partitions and
then apply the local rules of Section~\ref{sect:local} in the squares: for this, we order the squares (partially) by $(i,j)\preceq (i',j')$ if and only if $i\leq
i'$ and $j\leq j'$.

 \emph{We fix from now on a total order $\mathcal{O}$ extending this partial
order}. Every square $sq\neq (1,1)$ has then a predecessor $pred(sq)$, and every
square $sq\neq(\ell,\ell)$ has a successor $succ(sq)$. For $dir=\pm 1$, we define also $next(sq,dir)$ to be $succ(sq)$ if $dir=1$,  $pred(sq)$ if $dir=-1$), and to return ``undefined'' when $pred$ or $succ$ is not defined.

 Given $dir\in\{+1,-1\}$ and a square $sq$ of $G_\ell$:
\begin{itemize}
\item if $dir=1$, and $\mu,\nu,\lambda,C$ label $sq$ as in
 the definition of direct rules, we apply the corresponding rule;
\item if $dir=-1$, and $\mu,\nu,\xi$ label $sq$ as in the
definition of inverse rules, we apply the corresponding rule.
\end{itemize}

Let us call this procedure Apply\_local\_rule, and write
 loc:=Apply\_local\_rule(dir,$sq$) for the local rule that applies.

\medskip

\begin{figure}[!ht]
\centering
\includegraphics[width=0.9\textwidth]{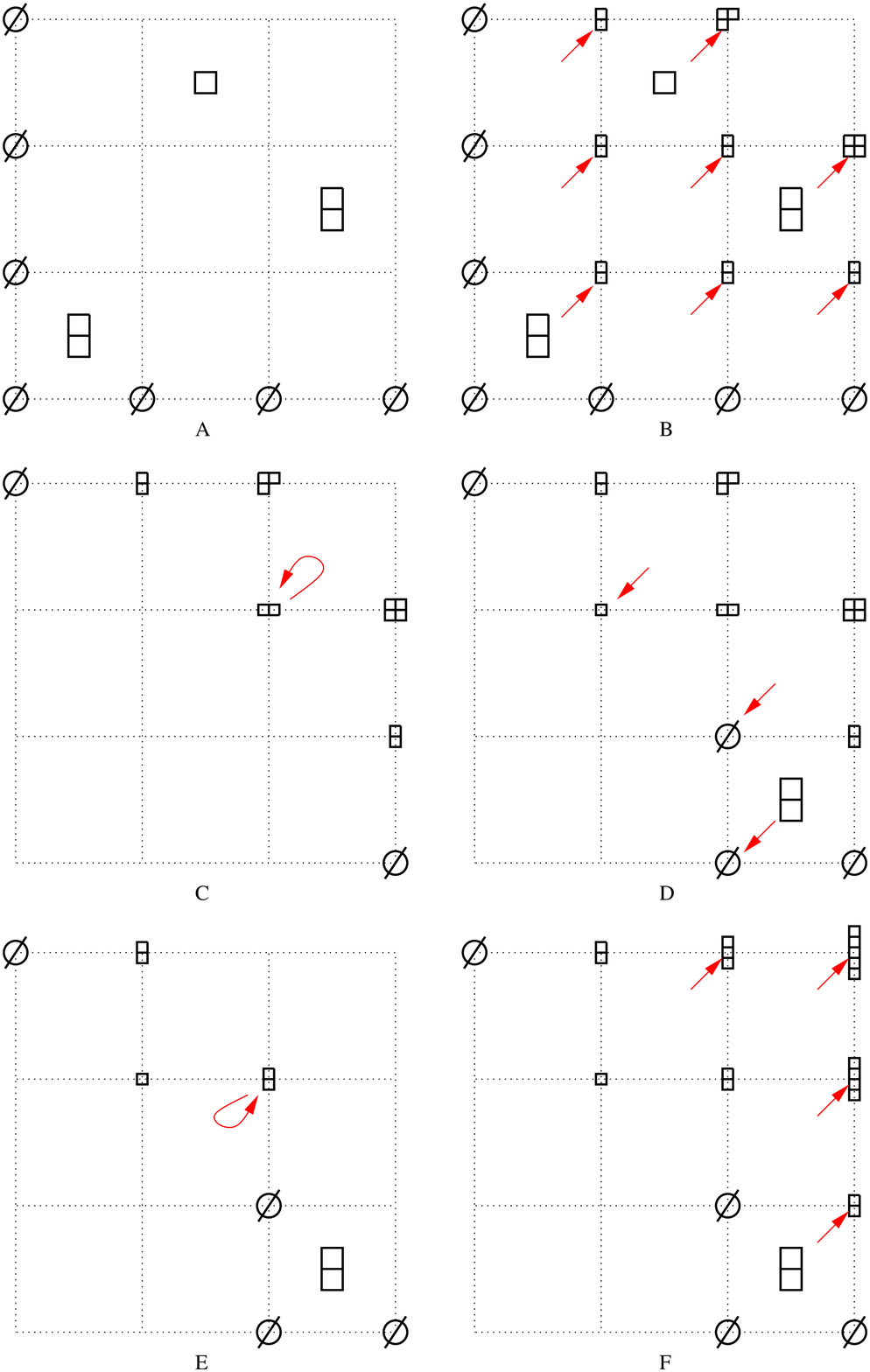}
 \caption{ \label{fig:bij} The algorithm $\phi$.}
\end{figure}

Now we want to go from local rules to a correspondence on the entire grid. Let us be given a hook permutation drawn $G_\ell$, in which we also label by $\emptyset$ (the empty partition) all vertices of $G_\ell$ on the bottom and left sides (see Configuration A on Figure~\ref{fig:bij}). We may now
describe the  bijection $\phi$ of Theorem~\ref{th:bij_perm}, which we do in an algorithmic fashion:\medskip
\medskip

 \fbox{
\begin{minipage}{0.9\textwidth}
   \noindent{\bf Algorithm $\phi:$}\\
Input: a hook permutation ($H,\sigma$).\\
Output: A pair ($P,Q$) of ribbon tableaux of the same shape.\\
\textbf{Begin}\\
\hspace{1cm}$sq:=(1,1)$; $dir:=1$;\\
 ~~{\bf repeat} \\
 ~~$loc$:=Apply\_local\_rule(dir,$sq$);\\
 ~~\textbf{If} ($loc\in$\{S,T\}) \textbf{then} $dir:=-dir$; \textbf{end if};\\
 $sq$:=next($sq$,dir);\\
 ~~\textbf{until} (sq=``undefined'');\\
\textbf{End}

\end{minipage}
}
\medskip

We will show that this algorithm is well defined and does not loop indefinitely:
it ends when $succ(\ell,\ell)$ is not defined, in which case the vertices on the
top and right side of $G_\ell$ are labeled by partitions forming two ribbon
tableaux of the same shape.
\medskip

{\noindent \bf Example:} we illustrate the algorithm on the example of Figure~\ref{fig:bij}. We choose the total order $(i,j)<(i',j')$ if $j<j'$, or
$j=j'$ and $i<i'$. A square $c'$ is thus larger than a square $c$ if it is above
in the same column, or if $c'$ is in a column to the right of $c$.

We start with the hook permutation given on $G_3$ (configuration A). We apply direct local
rules to reach configuration B, the rules being successively
D2, D4, D4, D4, D1, D2, D4 and D2. Now the rule S applies, and the direction changes
(configuration C). Note that we deleted all contents of visited squares as well
as the label of their bottom left corner; we did not show this for direct rules
in order to keep the number of pictures reasonable. From this configuration C, we apply inverse rules I5,I2 and I5 successively to
reach configuration D. There rule T applies, and the direction changes. We
finally reach configuration F with the rules D3,D2,D3 and D6. We can finally
read off the ribbon tableaux $P$ and $Q$ respectively on the right and top sides
of the grid.
\medskip

We will prove the correction of the bijection in Appendix~\ref{sect:appa}; we can already give the proof of its corollary: \medskip

{\bf Proof of Corollary~\ref{corol:enum_pairs}}: For the first formula, the signed bijection of Theorem \ref{th:bij_perm} implies that the left hand side is equal to the number of hook permutations of type $(c,c')$; this number is zero
unless $c$ can be obtained from  $c'$ by permuting some of its parts. When $\widetilde{c}=\widetilde{c'}$, we have to compute, for a part $i$ appearing
$j_i$ times, how many hook permutations of length $j_i$ exist so that all hooks
have size $i$; this number is clearly $i^{j_i}j_i!$. We obtain finally the right
hand side by multiplying such terms for all sizes $i$. The second equality of Corollary \ref{corol:enum_pairs} is an immediate consequence of Theorem
~\ref{th:bij_perm} and of Proposition~\ref{prop:enumhookperm}.

\subsection{Bijection for hook involutions}

 The Schensted correspondence is known to restrict to a bijection between involutions and standard tableaux, and we wish to extend this to hook involutions and ribbon tableaux. But whereas this restriction is immediate with Fomin's version of the Schensted correspondence, some extra work has to be done here: first, the bijection defined by the algorithmic procedure depends on the fixed total order specified on the squares of the grid, so the procedure can not be done symmetrically in general. Then, notice that hook involutions have signs, whereas hook permutations are defined as positive: we must prove that these signs are preserved by the bijections. We will deal with these two extra difficulties and give the proof of Theorem~\ref{th:bij_invol} in Appendix~\ref{sect:appa}.
\medskip

To end this section, we formulate in terms of signed sets the properties which are crucial in the proof of Theorem \ref{th:bij_invol}. Define $\U_i(\mu)$ as the set $\U_i(\mu,\mu)$
but with sign $\eps(\lambda/\mu)$ for $\lambda\in \U_i(\mu)$;
 similarly, $\D_i(\mu)$ is $\D_i(\mu,\mu)$ with sign
$\eps(\mu/\lambda)$.

\begin{proposition}[Shimozono and White~\cite{ShiWhi}] \label{prop:shiwhi}
Consider $\lb 0,i-1\rb$ as a signed set with $sgn(h)=(-1)^h$. Then 
 there is a sign preserving bijection between $\U_i(\mu)$ and
$\D_i(\mu)\sqcup \lb 0,i-1\rb$.
\end{proposition}

\subsection{Relation with previous works}

As already mentioned, Theorem~\ref{th:bij_perm}
was partly demonstrated by White~\cite{WhiteBij} (in this article,
 White evokes the possibility of extending his ideas to obtain
the form in which we gave it). In~\cite{StanWhi}, the authors
 notice that if one only considers hook permutations with
 all hooks of size $k$, and ribbon tableaux with all ribbons
 of size $k$, then rules S and T can never be applied: so in
 this case we have a (sign-preserving) bijection between ribbon tableaux
 and $k$-colored permutations.

As a matter of fact, rules D2 and D3 are not used in these two articles, but
alternative rules that do \emph{not} preserve signs in the sense of Proposition
\ref{prop:shiwhi}. Theorem~\ref{th:bij_invol} cannot thus be a consequence of
White's original work, but is based on the work of Shimozono and White
\cite{ShiWhi} (for ribbons and hooks of fixed size), in which the authors
introduce the operations $prev,next,first$ that are used to define rules $D2,D3$
and $I2,I3$.

\section{Algebraic approach} \label{sect:alg}

The previous section generalized the bijection of Robinson-Schensted;
 we will give an algebraic proof of the enumerative counterpart of
these results, that is Corollary~\ref{corol:enum_pairs}, in the spirit of Stanley~\cite{StanleyDiff} and Fomin
\cite{FominDual}.
\medskip

Let $\K$ be a field of characteristic zero. We consider
$\K\Y=\oplus_n\K\Y_n$, the vector space of all linear
combinations of partitions with coefficients in $\K$.
For $i$ a positive integer we define two linear operators
 $U_i$ and $D_i$ via their action on the basis of partitions:

\begin{definition} For $\la\in \Y$,
$$U_i\lambda=\sum_{r=\mu/\lambda\in Rib_i} \eps(r) \mu~~;
~~D_i\lambda=\sum_{r=\lambda/\mu\in Rib_i} \eps(r) \mu .$$
\end{definition}

$U_i$ and $D_i$ are endomorphisms of $\K\Y$, that send $\K\Y_{n}$ in
 $\K\Y_{n+i}$ and $\K\Y_{n-i}$
respectively. We note that these two operators were already defined by
Stanley~\cite{StanleyDiff} but in a different perspective.

 For $\la,\mu$ two partitions, we set $\langle \lambda,\mu\rangle =1$ if $\la=\mu$ and $0$ otherwise. We may then extend $\langle .,.\rangle $ to
$\K\Y\times\K\Y$
 by bilinearity, and we notice that $U_i$
 and $D_i$ are dual endomorphisms for this bilinear form; indeed
 $\langle U_i\lambda,\mu\rangle =\langle \lambda,D_i\mu\rangle $ since each
member in the
 equality is equal to $\eps(r)$ if $\la\subseteq \mu$ and
 $r=\mu/\la$ is a ribbon, and to 0 otherwise.

The fundamental relation between these endomorphisms is the following
 (we note $AB=A\circ B$ for composition):

\begin{proposition} 
\label{prop:relcoms}
For nonnegative integers $i,j$, we have
\begin{align} 
\label{eq:relcom1} D_iU_i=&U_iD_i+~ i\cdot\operatorname{Id}\\
\label{eq:relcom2} D_iU_j=&U_jD_i\quad\text{~~if~}i\neq j
 \end{align}
\end{proposition}

{\noindent{\bf Proof:~}} The first equality can be rewritten
\begin{eqnarray*}\langle D_iU_i\mu,\nu\rangle &=\langle U_iD_i+
i\cdot\operatorname{Id}\mu,\nu\rangle \\
\text{~~i.e.~}\langle D_iU_i\mu,\nu\rangle &=\langle U_iD_i\mu,\nu\rangle +
i\cdot\delta_{\mu,\nu}\\
\text{~~i.e.~}\langle U_i\mu,U_i\nu\rangle &=\langle D_i\mu,D_i\nu\rangle +
i\cdot\delta_{\mu,\nu}
\end{eqnarray*}

for $\mu$ and $\nu$ any two partitions, while the second equality is equivalent
to
$$\langle D_iU_j\mu,\nu\rangle =\langle U_jD_i\mu,\nu\rangle .$$

Those two equalities can be rephrased as
$|\U_i(\mu,\mu)|=|\D_i(\mu,\mu)|+i$ for all $\mu$, and
$|\U_i(\mu,\nu)|_\pm=|\D_i(\mu,\nu)|_\pm$ for
 $\mu\neq \nu$. But this is exactly the enumerative signification of Proposition~\ref{th:loc_bij}, which ends the proof. \qee

\medskip

The first of these two relations is characteristic of {\em
$i$-differential posets}, defined by Stanley
\cite{StanleyDiff,StanleyVarDiff}, and which are the basis
 of~\cite{FominDual} (let us remark that our operators involve
 signs, which is not the case in the classical case.)\medskip.

Let us now specialize $\K=\mathbb{Q}((q))$, the field of formal Laurent series
 in  $q$ (We could in fact limit ourselves to formal power
 series $\mathbb{Q}[[q]]$, and in general work with
 commutative rings of characteristic zero).

Let us define $\K^\Y=\prod_n\K\Y_n$ the vector space
of functions from $\Y$ to $\K$, which we will write as infinite linear
combinations of partitions,
with coefficients in $\K$. It is then possible to extend
  $\langle .,.\rangle $ to $\K\Y\times \K^\Y$ without difficulty, and one checks
that
   $D_i$ and $U_i$ extend also to endomorphisms of $\K^\Y$. Let us now define
operators $\UU$ and $\DD$:

\begin{definition} \label{def:UUDD}
$$\UU=\sum_{i} q^iU_i {~~;~~}\DD=\sum_{j} q^jD_j$$
\end{definition}

$\UU$ and $\DD$ are themselves endomorphism of $\K^\Y$. $\UU(\la)$ (respectively
$\DD(\la)$) is by definition the sum of all partitions $\mu$ such that $\mu/\la$
(resp. $\la/\mu$) is a ribbon, with a coefficient  $\eps q^k$ for a ribbon of
size $k$ and sign $\eps$ (note that $\UU$ does not have its image in $\K\Y$) .
This interpretation gives us immediately the following formula that explains our
interest in $\UU$ and $\DD$:

 \begin{equation}\label{eq1}
\sum_{\substack{\lambda\in \Y_n\\ P,Q\in RT_{\lambda,\ell}}}
 \varepsilon(P)\varepsilon(Q) =
   [q^{2n}]\langle \emptyset,\DD^\ell\UU^\ell\emptyset\rangle 
 \end{equation}

 So we want to calculate the right hand side of this equality;
 note first that Proposition~\ref{prop:relcoms} can be summarized
 in a single relation involving  $\UU$ and $\DD$:

\begin{equation}\label{eq:bigrel}
\DD\UU=\UU\DD+ \frac{q^2}{(1-q^2)^2}\operatorname{Id}
\end{equation}

\noindent{\bf Proof:~} One just has to take the definition
of $\UU$ and $\DD$, develop, and then use the relations 
of Proposition~\ref{prop:relcoms}.
 The coefficient of $\operatorname{Id}$ appearing is then seen
 to be $\sum_{i\geq 1} iq^{2i}$, which is equal to the given
 rational fraction. \qee\medskip

In order to compute the coefficient in $q$ in Equation~\eqref{eq1},
we need the following result of Stanley :

\begin{theorem}[Stanley~\cite{StanleyDiff}] \label{th:diffpos}
Suppose that two endomorphisms of a vector space $E$ verify
$DU=UD+rI$. Then for all positive integers $\ell$ we have
 \begin{equation}\label{eq:stan}
D^\ell U^\ell=(UD+rI)(UD+2rI)\cdots(UD+\ell rI)
\end{equation}
As a consequence, if $\widehat{O}$ is an element of $E$
such that $D\widehat{O}=0$, we have $\langle \widehat{O},D^\ell
\mathrm{}_{}U^\ell\widehat{O}\rangle =r^\ell\ell!$
\end{theorem}

 We have such a relation for $\UU$ and $\DD$ with
 $\widehat{O}=\emptyset$ and $r=q^2/(1-q^2)^2$ : it's just equation
\eqref{eq:bigrel}.
So the second formula of~\ref{corol:enum_pairs} is the
 consequence of the following computation:

$$r^\ell\ell!=\ell!\cdot q^{2\ell}\cdot\frac{1}{(1-q^2)^
{2\ell}}=\sum_{n\geq 2\ell}
\left[\binom{n+\ell-1}{2\ell-1}\cdot \ell!\right] q^{2n}.$$

\medskip

Let us now turn to the first formula of Corollary
\ref{corol:enum_pairs}.  $U_{\mu_\ell}\cdots
U_{\mu_1}\emptyset$ is the linear combination of ribbon tableaux
 of content $\mu$, with the sign of the tableau as coefficient.
 So the left hand side of
\ref{corol:enum_pairs} is $\langle U_{\nu_\ell}\cdots
U_{\nu_1}\emptyset,U_{\mu_\ell}\cdots U_{\mu_1}\emptyset\rangle $, which by
duality is
 $\langle \emptyset,D_{\nu_1}\cdots
D_{\nu_\ell}U_{\mu_\ell}\cdots U_{\mu_1}\emptyset\rangle $.

\begin{lemma} Let $\mu,\nu$ be two compositions with $\ell$ parts.
Then $$\langle \emptyset,D_{\nu_1}\cdots D_{\nu_\ell}U_{\mu_\ell}\cdots
U_{\mu_1}\emptyset\rangle  = \nu_{\ell}\times\sum_\rho
\langle \emptyset,D_{\nu_1}\cdots D_{\nu_{\ell-1}}U_\rho\emptyset\rangle ,$$
where $\rho$ goes through the multiset of all compositions of
length $\ell-1$ obtained by deleting a part of $\mu$ of size $\nu_{\ell}$.
\end{lemma}

\noindent{\bf Proof of the lemma:~} It is just an iterated use of Proposition \ref{prop:relcoms}, with $D_{\nu_\ell}$ and $U_{\mu_i}$, for $i$ equals $\mu_\ell$,\ldots,$\mu_1$ successively. Two cases can occur : if $\mu_i\neq\nu_\ell$, then $D_{\nu_\ell}$ and $U_{\mu_i}$ commute; otherwise, for
each index $i$ such that $\mu_i=\nu_\ell$, a term $\nu_{\ell}\langle
\emptyset,D_{\nu_1}\cdots D_{\nu_{\ell-1}}U_\rho\emptyset\rangle $ has to be
added, where $\rho$ is the composition with the part $\mu_i$ deleted in $\mu$.
Then when $D_{\nu_\ell}$  $U_{\mu_i}$, the scalar product is zero since
$D_{\nu_\ell}\emptyset=0$. \qee\medskip

The proof of the first part of the Corollary~\ref{corol:enum_pairs} is now a simple induction on $\ell$ based on the preceding Lemma.


\section{Columns of the character table of $\St_n$}
 \label{sect:char}

Ribbon tableaux are known to have a strong connection with the representation
theory of the symmetric group, which we will now recall briefly. Let $n$ be a
positive integer, $\la$ a partition of size $n$; we note $\chi^\la$ 
the irreducible character of $\St_n$ indexed by $\la$ (for
more information on these topics, see for instance~\cite{FultonBook,SaganBook}
that have a combinatorial approach). Let also $\chi_\mu^\la$ be the value of
this character on a permutation of cycle type $\mu$: this means that the permutation
has $m_i$ cycles of length $i$ for each $i$ if $\mu=(1^{m_1}2^{m_2}\cdots)$. The
\emph{Murnaghan-Nakayama rule}~\cite{Mur,Nak} states that:

\begin{theorem}  Let $\mu,\nu$ be two partitions of the same size $n$.
Then

$$\chi_\mu^\lambda=|RT_{\la,\mu}|_\pm$$

\end{theorem}

This rule gives a combinatorial interpretation of $\chi_\mu^\lambda$,
 and shows in particular that it is an integer. We will now show that
Theorem~\ref{th:bij_invol} is adapted to study the {\em column sums}
of the character table.


\subsection{A formula for $\sum_\lambda \chi_\mu^\lambda$}

Define $C(\mu)=\sum_\lambda \chi_\mu^\lambda$ , the sum of all entries of column $\mu$ in the character table of $\St_n$. By the Murnaghan-Nakayama rule, $C(\mu)$ is equal to the signed sum of all ribbon tableaux of content $\mu$. By Theorem~\ref{th:bij_invol}, this last quantity is itself equal to the signed sum
of hook involutions of content $\mu$. The preceding result can thus be summed up by $C(\mu)=|RT_{\mu}|_{\pm}=|\HI|_{\pm}$. Using Corollary \ref{corol:inv_hook_inv}, we finally obtain
\begin{equation} \label{eq:cmu}
C(\mu)=|{\HI_{spec}(\mu)}|.
\end{equation}

This shows in particular that $C(\mu)$ is a {\em nonnegative } integer; the
following theorem gives a formula for the exact value of this integer.

\begin{theorem} \label{th:colsum}
Let $\mu=(1^{m_1}2^{m_2}\cdots)$  be a partition.Then
$C(\mu)=\prod_{i>0} c_{i,m_i}$ with:
\begin{eqnarray*} c_ {i,m_i}=& 0&\text{ if }i\text{ is even and
}m_i\text{ is odd;}\\
 =&(m_{i}-1)!!\cdot i^{m_{i}/2}&\text{ if }i\text{ is even and
 }m_i\text{ is even;}\\
     =& \sum_{k=0}^{\left\lfloor\frac{m_{i}}{2}\right\rfloor}
 \binom{m_{i}}{m_{i}-2k}\cdot (2k-1)!!\cdot i^k&\text{ if i is
 odd.}
 \end{eqnarray*}
\end{theorem}

We will exhibit two proofs: first a bijective one,
 and the second algebraic, using the tools of Section~\ref{sect:alg}.
\medskip

\noindent {\bf First Proof: } The computation of
$|{\HI_{spec}(\mu)}|$ for general $\mu$ reduces clearly to the
case $\mu=[i^{m_i}]$ where $\mu$ has only one part size.

In this case an element of ${\HI_{spec}(\mu)}$ is an
involution on $\lb 1,a_k\rb$ with a choice of a hook of size $a_k$
for each cycle of length $2$. Remembering that elements of
${\HI_{spec}(\mu)}$ have no fixed points corresponding to even
parts, we obtain easily the above expression for the
 coefficient $c_{i,m_i}$, and the proof is complete. \qee
\medskip

\noindent {\bf Second proof: } We now give a proof that does
not use Theorem~\ref{th:bij_invol} (and thus also not
 the Equality~\eqref{eq:cmu}). For this we need the following algebraic
consequence of Proposition~\ref{prop:shiwhi} in terms of the operators
$D_i$ and $U_i$ (considered as endomorphisms of $\K^\Y$ ); we
note $\Y$ the vector $\sum_{\lambda\in \Y} \lambda \in \K^\Y$, and
$o_i$ is 1 when $i$ is odd and 0 otherwise.

\begin{proposition} \label{prop:sigalg}
 For all $i\geq 1$, we have $ D_i\Y=U_i\Y+o_i\cdot \Y$.
\end{proposition}

\noindent {\bf Proof: } Take the scalar product of each member
of the equality with a partition  $\la$, and remembering that
 $U_i$ and $D_i$ are dual operators, the result is equivalent to
 $$\sum_{\mu\in\mathcal{U}(\la)}\eps(\mu/\la)  =
\sum_{\mu\in\mathcal{D}(\la)}\eps(\la/\mu)   +o_i.$$

This is an immediate corollary to Proposition~\ref{prop:shiwhi}.\qee \medskip

By the Murnaghan-Nakayama rule, we have $C(\mu)=\langle
D_\mu\Y,\widehat{\emptyset}\rangle $.
We will use the relations of Propositions~\ref{prop:relcoms}
and~\ref{prop:sigalg} to compute this scalar product.

\begin{lemma} We have the following formulas:
\begin{itemize}
\item[1.] For $m\geq 2$ and $i\geq 1$, $D_i^m\Y=o_i\cdot
D_i^{m-1}\Y+(m-1)i\cdot D_i^{m-2}\Y+U_iD_i^{m-1}\Y.$ 

\item[2.] For $m\geq 1$ and $i\geq 1$, $D_i^m\Y=c_{i,m}\Y+U_iA_{i,m}\Y$, where
$c_{i,m}$ is defined in Theorem~\ref{th:colsum}, and $A_{i,m}$ is an
endomorphism of $\K^\Y$. 
\item[3.] For $m\geq 1$ and $i\geq 1$, $\langle D_\mu D_i^m\Y,\emptyset\rangle =c_{i,m}\langle
D_\mu \Y,\emptyset\rangle $~if all parts of $\mu$ are greater than $i$.
\end{itemize}
\end{lemma}

\noindent {\bf Proof of the lemma:~}
 We have  $D_i^m\Y= o_i\cdot D_i^{m-1}\Y+D_i^{m-1}U_i\Y$ thanks to Proposition
~\ref{prop:sigalg}. Using $m-1$ times the relation
 $D_iU_i=U_iD_i+i\cdot I$, point 1. of the lemma is proved.

  By an immediate induction on 1., we can write for all $m\geq 2$ that
$D_i^m\Y=b_{i,m}\Y+U_iB\Y$ for a certain endomorphism $B$ and an integer
$b_{i,m}$ necessarily equal to $\langle D_i^m\Y,\emptyset\rangle $. Substituting
in 1., and taking the coefficient of $\emptyset$ in each member, we obtain
$b_{i,m}=o_ib_{i,m-1}+(m-1)ib_{i,m-2}$. The numbers $c_{i,m}$ verify the same
recurrence relation, as can be easily seen, directly or by the combinatorial interpretation 
given in the first proof. Since we have in addition
$b_{i,0}=c_{i,0}=1$ and $b_{i,1}=c_{i,1}=o_i$, it follows $b_{i,m}=c_{i,m}$ for
all $i,m$  and point 2. is proved .

 Finally, thanks to point 2., the left member of 3. is equal to $$c_{i,m}\langle
D_\mu\Y,\emptyset\rangle +\langle D_\mu U_iA_{i,m}\Y,\emptyset\rangle ,$$ Since $D_\mu$ commutes with $U_i$ by Equation~\eqref{eq:relcom2}, the second term is then $0$ because the image of $U_i$ has null intersection with $\K\emptyset$, and the lemma is proved.\qee

The proof of Theorem~\ref{th:colsum} is now immediate
 by induction on the number of different part sizes of 
 $\mu$, using formula 3. in the previous lemma.\qee

\subsection{Other enumerations of $C(\mu)$}
The formula of Theorem ~\ref{th:colsum} is not new, but the proof above is (to the best of our knowledge) the first fully bijective proof of it based on the Murnaghan Nakayama rule. Let us mention two other places in the literature where this result is shown, and show the equivalence to our formulation.

The computation of $C(\mu)$ is an exercise in Macdonald's book
\cite[p.122, ex.11]{Macdonald}, and relies on symmetric function techniques. It is proved that $C(\mu)$ is equal to the product $\prod_{i\geq 1} a_i^{(m_i)}$, where $a_i^{(m)}$ is the coefficient of  $t^m/(m!)$ in $\exp(t+\frac{1}{2}it^2)$ (\emph{respectively} $\exp(\frac{1}{2}it^2)$) if $i$ is odd ({\emph{resp.} even)}. Through an expansion of the series, one checks easily that $a_i^{(m)}$
is indeed equal to the coefficient $c_{i,m}$ of Theorem~\ref{th:colsum}.

Another proof can be found in Exercise $7.69$ of Stanley's book
\cite{StanleyEnum2}; the proof is based on a general result
in character theory, whose specialization to the symmetric group is the following theorem:

\begin{theorem}[\cite{Isaacs,StanleyEnum2}] \label{th:hitoroot}
Let $\s$ be a permutation of $\lb 1,n\rb$ with cycle type $\mu$.
Then $C(\mu)$ is equal to the number of square roots of $\s$ in $\St_n$, i.e. to the number of permutations $\tau\in \St_n$ such that
$\tau^2=\s$.
\end{theorem}

\noindent {\bf Proof: } Thanks to the formula~\eqref{eq:cmu}, we can prove this result by exhibiting a bijection $HiToRoot$ between $\HI_{spec}(\mu)$ and $\{\tau~|~\tau^2=\s\}$.

We consider each cycle of $\s$ as a word $x\s(x)\s^2(x)\ldots$ where $x$ is minimal in its orbit, and decompose $\s$ canonically in the form $[c_1^{(1)}\cdots c_1^{(m_1)}][c_2^{(1)}\cdots c_2^{(m_2)}]\cdots $, where $c_i^{(k)}$, $k=1\ldots m_k$, are the cycles of length $i$ written in increasing
order of their minimal elements. For instance, the permutation $574389216$ has
cycle type $\mu=(3,2^3)$ and will be written $[(27)(34)(69)][(158)]$.

\begin{lemma}
Let $c_1,c_2$ be two disjoint cycles of length $m$ in $\St_n$.\\
$\bullet$ There exist exactly $m$ cycles $c$ of length $2m$ in $\St_n$
such that $c^2=c_1c_2$.\\
$\bullet$ If $m$ is odd, there exist a unique cycle $c$ of length
$m$ in $\St_n$ such that $c^2=c_1$.
\end{lemma}

The proof of this lemma is immediate. The $m$ cycles of the first part of the
Lemma will be denoted $root(c_1,c_2,j)$, $j=0\ldots m-1,$ and the unique cycle
of length $j$ determined by the second part is $root(c_1)$.

 Let us now fix an element $I\in {\HI_{spec}(\mu)}$; it is equivalent to a
sequence of hook involutions $I_j, j=1\ldots k$ where $I_j$ is element of
${\HI_{spec}(j^{m_j})}$. Write $t_j$ ({\em resp. } $f_j$) for the number of
transpositions ({\em resp.} fixed points) in the involution $I_j$.  Let us also
write $(x_s,y_s), s=1\ldots t_j$ these transpositions  and $z_t, t=1\ldots f_j$ for
the fixed points. Finally let $h_s\in\lb 0,j-1\rb$ be the height of the hook
associated to the transposition $(x_s,y_s)$.

 We now associate to $I_j$ the $t_j$ cycles of length $2j$ defined by 
$root(c_j^{x_s},c_j^{y_s},h_s),s=1\ldots t_j,$ as well as the cycles of length
$j$ $root(c_j^{(z_t)}),t=1\ldots f_j.$ The product of all these disjoint cycles
for all indices $j$ form a permutation, which is the desired root $HiToRoot(I)$.\qee
\medskip

We can also use the formula to answer the question : for a given
integer $k$, what are the partitions $\mu$ such that the column
sum $C(\mu)$ is equal to $k$ ? Let $\mathcal{O}\mathcal{D}$ be the
set of partitions with odd distinct parts. The answers for the
first integers are:

\begin{itemize}

\item $C(\mu)=0$ if and only if $\mu$ has at least an even part with odd
multiplicity;

\item $C(\mu)=1$ if and only if $\mu\in \mathcal{O}\mathcal{D}$.

\item $C(\mu)=2$ if and only if $1$ has multiplicity $2$ and $\mu-1^2\in
\mathcal{O}\mathcal{D}$, or $2$ has multiplicity $2$ and
$\mu-2^2\in\mathcal{O}\mathcal{D}$

\item $C(\mu)=3$ has no solution.

\item $C(\mu)=4$ if and only if $3$ has multiplicity $2$ and $\mu-3^2\in
\mathcal{O}\mathcal{D}$, or $4$ has multiplicity $1$ and
$\mu-4^1\in \mathcal{O}\mathcal{D}$, or $2$ and $1$ have
multiplicity $2$ and $\mu-1^22^2\in \mathcal{O}\mathcal{D}$.

\end{itemize}

The number of solutions to $C(\mu)=0$ is sequence A085642 in
Sloane's Online Encyclopedia~\cite{Online}. The article
\cite{BesOls} proves that another family of partitions is in
bijection with $\mathcal{O}\mathcal{D}$, namely the partitions
with at least one part congruent to $2$ modulo $4$.

\section{Extensions} \label{sect:ext}

In this last section we sketch three different directions for which the ideas of this work can be applied.

\subsection{Combinatorial proof that characters are class functions}

Stanton and White~\cite{StanWhi} show combinatorially that, if $c$
and $c'$ are $2$ compositions verifying $\widetilde{c}=\widetilde{c'}$, 
then
\begin{eqnarray}
\label{eq:invar}
|RT_{\lambda,c}|_\pm=|RT_{\lambda,c'}|_\pm 
\end{eqnarray}
 
 This expresses in fact that the value of the character $\chi^\lambda$ on a
permutation $\s$ depends only on the conjugacy class of $\s$: it is a {\em class
function}. In this Section we give local rules that realize Stanton and White's result,
building on Fomin's version of jeu de taquin explained for instance in his
appendix to Stanley's book~\cite{StanleyEnum2}. We have the following
proposition, whose proof can be easily done by using the encoding of partitions
by infinite sequences explained in Appendix~\ref{sect:appb}.

\begin{proposition}
\label{prop:last}
Let $\la,\mu,\xi$ be three partitions, such that $\mu/\la$ and $\xi/\mu$ are 
nonempty ribbons.  Then exactly one of the two following cases occur:
 \begin{enumerate}
\item either there exists  $\nu$ such that $\nu/\la$ and $\xi/\nu$ are ribbons
of respective sizes  $|\nu/\la|=|\xi/\mu|$ and $|\xi/\nu|=|\mu/\la|$, or
\item there exists $\widehat{\mu}$ such that $\widehat{\mu}/\la$ and
$\xi/\widehat{\mu}$ are ribbons of respective size $|\widehat{\mu}/\la|=|\mu/\la|$
and $|\xi/\widehat{\mu}|=|\xi/\mu|$. 
\end{enumerate}

Furthermore, we have the sign relations
$\eps(\nu/\la)\eps(\xi/\nu)=\eps(\mu/\la)\eps(\xi/\mu)$ in the first case,
and $\eps(\widehat{\mu}/\la)\eps(\xi/\widehat{\mu})=-\eps(\mu/\la)\eps(\xi/\mu)$
in the second case.
\end{proposition}

\begin{figure}
      \psfrag{la}{$\la$}
      \psfrag{mu}{$\mu$}
      \psfrag{xi}{$\nu$}
      \psfrag{nu}{$\xi$}
 \psfrag{11}{$(1,1)$}\psfrag{12}{$(1,2)$}\psfrag{22}{$(2,2)$}\psfrag{13}{$(1,3)$}
 \psfrag{33}{$(3,3)$}\psfrag{23}{$(2,3)$}
 \includegraphics[height=5cm]{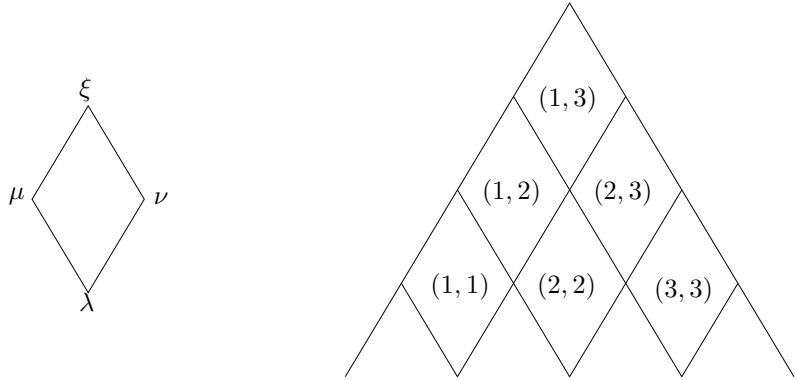}
\caption{Lozenge for local rules and associated grid for the global correspondence.
\label{fig:ruleevac}}
\end{figure}

We can now define the local rules: given $\la,\mu,\xi$ as in the
proposition, we draw them on a lozenge as on the left of Figure
\ref{fig:ruleevac}. If the first case of the proposition occurs, we erase $\mu$,
and write $\nu$ on the right vertex of the lozenge; otherwise we replace $\mu$
by $\widehat{\mu}$. We also define inverse local rules by a simple vertical symmetry.
Finally, we define the {\em trivial} local rule which consists of simply moving the 
partition $\mu$ from one side to the other.

Now we go from the local to the global as we did in Section~\ref{sect:bij}.
Consider the grid $P_\ell$ on the right of Figure
\ref{fig:ruleevac}, made of lozenges. We attach to each lozenge coordinates $(i,j)$, with 
$1\leq i\leq \ell-1$ and $i\leq j\leq \ell-1$ in the manner shown in the example.
We fix a total order on lozenges, such that each lozenge $(i,j)$
 has to be bigger than the two lozenges on its top left and bottom left, namely
 $(i-1,j)$ and $(i,j-1)$ when they are defined. For the examples, we will use
 the following linear order: $(i,j)$ is greater than $(i',j')$ if $i>i'$, or $i=i'$ and $j>j'$. 

Let us fix two compositions $c,c'$ of length $\ell$, such that
$\widetilde{c}=\widetilde{c'}$. We will represent elements of $RT_{\lambda,c}$
on the left side of $P_\ell$, as chains of partitions labeling the vertices from bottom to
top, and elements of $RT_{\lambda,c'}$ on the right side in the same fashion.

Now we want to select a subset of the lozenges so that when local rules are applied in the grid, we
obtain indeed a correspondence between $RT_{\lambda,c}$ and $RT_{\lambda,c'}$.
Fix $\sigma$ a permutation such that $\sigma(c)=c'$, and mark certain lozenges of 
$P_\ell$, in the following 
way: for $1<j\le l$, if $\sigma_j=k$,
then mark the lozenges $(i,j-1)$, $i\le l-k$.

\begin{figure}
\psfrag{1}{$1$}\psfrag{2}{$2$}\psfrag{3}{$3$}\psfrag{4}{$4$}
 \includegraphics[height=4cm]{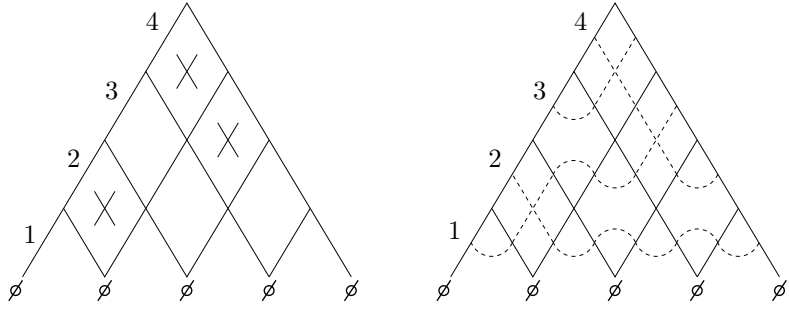}
\caption{Marking of the lozenges for the permutation $3142$, and resulting
effect on the size of ribbons through the correspondence.\label{fig:marks}}
\end{figure}

For instance, consider $c=(1,3,2,1)$ and $c'=(3,1,1,2)$, and fix the permutation
$w=3142$ which indeed verifies $w(c)=c'$; the corresponding marking is represented on the
left hand side of Figure~\ref{fig:marks}.

\begin{figure}
 \includegraphics[height=4cm]{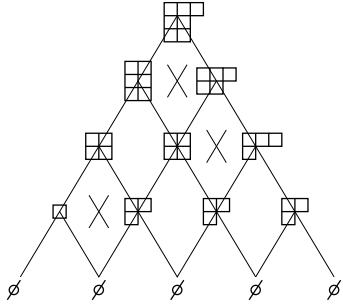}
\caption{Example of the correspondence in the proof of Equation~\eqref{eq:invar}.\label{fig:invar}}
\end{figure}

Now we start from a ribbon tableau represented on the left hand side, and go on
performing local rules as in the beginning of Section~\ref{sect:bij}: the rule
being that we perform nontrivial local rules in the marked lozenges, while in
the other (unmarked) ones we will use the trivial local rule. An example is
given on Figure~\ref{fig:invar}; the definition of the marked lozenges in the
previous paragraph is thus made so that the size of ribbons ``match'' between $c$ and
$c'$, thanks to Proposition~\ref{prop:last}; this can be visualized on the right hand side of Figure~\ref{fig:marks}.

 We remark that in the case when all lozenges are marked with a cross, then this corresponds to a generalization of Sch\"utzenberger's involution which is the case where all ribbons are of size $1$.

\subsection{Other correspondences based on the graph $GR$}

In the correspondence of Theorem~\ref{th:bij_perm}, we considered pairs of ribbon tableaux of the same shape. As already noticed, these are very special paths in the ribbon graph GR: they start and end at $\emptyset$, going up $\ell$ steps and then down $\ell$ steps. The same ideas work to build correspondences for other kinds of paths, and we give an example that is well known in the standard case.

 We consider the paths in $GR$ of length $2\ell$, that start and end at $\emptyset$, and which
possess $\ell$ steps up and $\ell$ steps down (but in no imposed order). These
are called oscillating tableaux (of shape $\emptyset$) in the case where all
ribbons are of size $1$, and we will thus call these paths \emph{oscillating
ribbon tableaux}.

\psfrag{sub}{$\subset$}   \psfrag{sup}{$\supset$}
\begin{center}
 \includegraphics[height=1cm]{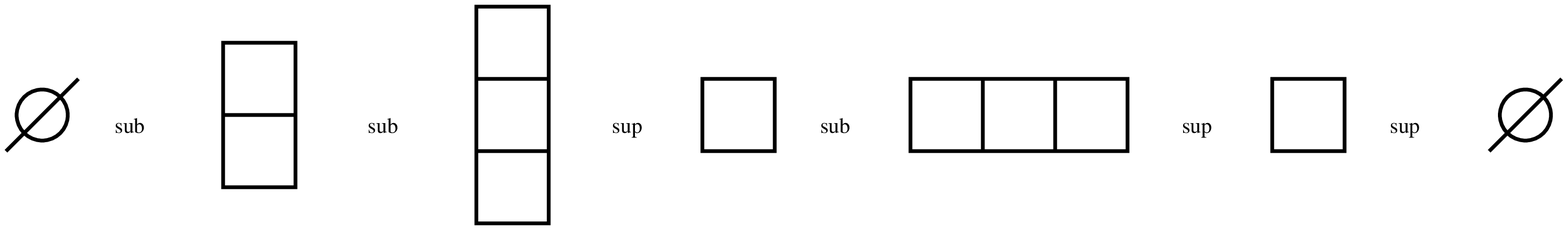}
\end{center}

The {\em size} of such a path is the \emph{half sum} of the
sizes of the $2\ell$ edges (i.e. ribbons), and the sign is the product of the
signs of those ribbons. The oscillating tableau above has 
size $5$, sign $+1$ and length $6$. Let us note $Osc_{n,\ell}$ this signed set; we can 
then prove the following formula:
$$|Osc_{n,\ell}|_\pm=(2\ell-1)!!\binom{n+\ell-1}{\ell-1}$$

This can be done in two ways, algebraic and bijective, following the steps of
what was done in the case of pairs of ribbon tableaux of the same shape.
 \medskip
 
 The algebraic way to prove the identity is to notice that the quantity
$|Osc{n,\ell}|_\pm$ can be expressed as the coefficient of $q^{2n}$ in the
series $\langle \emptyset,(\DD+\UU)^{2\ell}\emptyset\rangle $. Now this series
is $(2\ell-1)!!(\sum_iiq^{2i})^\ell$ (this a consequence of Corollary 2.6 (a) of
\cite{StanleyDiff}); then the end of the proof goes as for the second equality of
Corollary~\ref{corol:enum_pairs}.
\medskip

The bijective way consists in constructing a signed bijection between
$Osc_{n,\ell}$ and {\em hook matchings} of $\lb 1,2\ell\rb$ with size $n$:
these are perfect matchings on $\lb 1,2\ell\rb$ such that to each pair
 $\{i,j\}$ of the matching is associated a hook $H_{\{i,j\}}$, such that the sum
of the sizes of the $\ell$ hooks
is $n$. Since there are $(2\ell-1)!!$ perfect matchings, Proposition
\ref{prop:enumhookperm} shows that there are 
$(2\ell-1)!!\binom{n+\ell-1}{\ell-1}$ such hook perfect matchings.

\begin{figure}
\begin{center}
 \includegraphics[height=8cm]{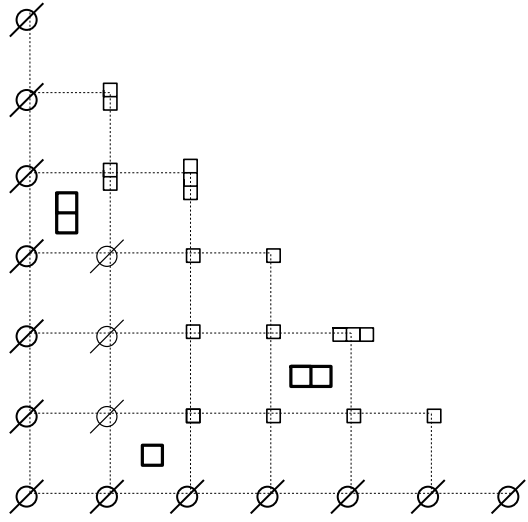}
\end{center}
\caption{Correspondence for oscillating ribbon tableaux \label{fig:oscil}}
\end{figure}

 The bijective correspondence between $Osc_{n,\ell}$ and the hook matchings
is done as in Roby~\cite{Roby}. We illustrate this on an example on Figure
\ref{fig:oscil}. Instead of the grid $G_\ell$, we perform the bijection on a
grid $T_\ell$ illustrated on the Figure by dashed lines for $\ell=3$.

 Hook matchings can be represented by labeling by $\emptyset$ the bottom and
left side, and for each pair $\{i,j\}$ of the matching, the corresponding hook
is drawn in the square of column $i$ from the left and row $j$ from the top. In
the example, the matching is then $\{\{1,3\},\{2,6\},\{4,5\}\}$. Oscillating
ribbon tableaux
$(\la_0=\emptyset,\la_1,\ldots,\la_{2\ell-1},\la_{2\ell}=\emptyset$ are
represented on the outside corners of the north east border, from top left to
bottom right ( moreover, in each of the corresponding inside corners, one draws
the smallest shape between $\la_i$ and $\la_{i+1}$).

Now a signed correspondence goes along the exact same lines as what we did in
Section~\ref{sect:bij} for pairs of ribbon tableaux of the same shape: we fix a total
order on the squares of $T_\ell$, and apply local rules in these squares,
changing directions when we encounter a rule $S$ or $T$. The example of Figure
\ref{fig:oscil} does not present any occurrence of those last rules for the sake
of simplicity.

\subsection{Layered graphs in duality}

The techniques used in this article are a generalization of Fomin's framework developed in \cite{FominDual,FominSchen}; we now present our theoretical setting in this paragraph.

 Consider a graph $G=(V,E)$ with a sign function on the edges
 $\eps:E\rightarrow\{+1,-1\}$. Suppose that $V$ is the disjoint union
  of finite sets $V_i$, $i\in \mathbb{N}$  where $V_0$ is a singleton $\{O\}$.  We will say that $G$ is  {\em layered} graph (with zero).

  Let now $U_i,D_j$ be the endomorphisms of $\K V$ defined for $v\in V_k$ by $U_i(v)=\sum_e \eps(e)v'$ where $e$ goes through all edges from $v$ to $v'\in V_{k+i}$, and $D_j$ is defined dually from $\K V_k$ to $\K V_{k-i}$.\smallskip

Let us say that the layered graph $G$ is {\em self dual} if
there exist nonnegative integers $\alpha_i$ such that:
\begin{align*}
D_iU_i=&U_iD_i+~ \alpha_i\cdot\operatorname{Id}\\
D_iU_j=&U_jD_i\quad\text{~~if~}i\neq j
\end{align*}

Fomin's framework of self dual graded graphs is the
the case where edges exist only between consecutive
levels $V_i$ and $V_{i+1}$, and when the sign function is constant equal to $1$.

It is then possible to use the algebraic techniques of this work to study the enumeration of paths in such graphs. The relations above correspond to certain equalities of signed cardinals, as in the graph $GR$. If signed bijections proving these equalities are fixed, then we can determine global correspondences in the same way. But it obviously remains to see if there exists interesting examples to which this theory can be applied.

An important remark is that this generalizes only Theorem
\ref{th:bij_perm} and its Corollary; to give an extension of Theorem
\ref{th:bij_invol}, one needs additional local properties on the self dual layered graph, analogous to Proposition~\ref{prop:sigalg} in the case of $GR$.

\appendix

\section{Partitions as sequences of zeros and ones.}
\label{sect:appb}

In this appendix we will show how the encoding of partitions by words is well suited to the study of the operations on ribbons of Section \ref{subsect:oper}. We follow van Leeuwen~\cite{VanLeeuEdge} for notations. Let $\la$ be a partition, and $\de(\la)\in\{0,1\}^\Zm$ the sequence defined by the following procedure: we extend the top and left borders of a Ferrers diagram to infinity, and read the lower right boundary from bottom to top, recording $1$ for every vertical edge encountered, and $0$ for the horizontal edges.

\begin{figure}
\centering
\includegraphics[width=.4\textwidth]{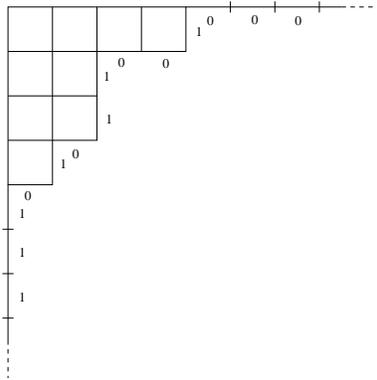}
\caption{\label{Prelim:fig:edgeseq} Infinite word encoding a partition.}
\end{figure}

For instance, the partition $(4,2,2,1)$ has for coding word
 $(\cdots 1110101|1001000\cdots ),$ cf. Figure~\ref{Prelim:fig:edgeseq};
the sign ``$|$'' separates the parts of the border below and above the
diagonal of the diagram, and we consider that nonnegative indices of
$\delta(\lambda)$ are those on the right of $|$. Notice that the encoding
 sequences have the following characteristic properties (see
\cite{VanLeeuEdge}):

\begin{enumerate}
\item  they differ from $(\cdots 1111|0000\cdots )$
(corresponding to the empty partition) in a finite number of positions;
 \item the number of $0$s to the left of $|$
is equal to the number of $1$s to its right.
\end{enumerate}

Now, ribbons addable to $\la$ (respectively removable from $\la$)
 are in bijection with pairs of indices $(i,j)$ in $\Zm$ where $i<j$
such that $\de_i(\la)=1$ and $\de_j(\la)=0$ (resp. $\de_i(\la)=0$ and
$\de_j(\la)=1$): $i$ indicates the position of the head, and
 $j$ the position of the tail. The partition $\mu$ obtained
 after addition or removal of the ribbon is the result of the
 exchange of $0$ and $1$ at positions $i$
and $j$.

\begin{figure}
\includegraphics[width=\textwidth]{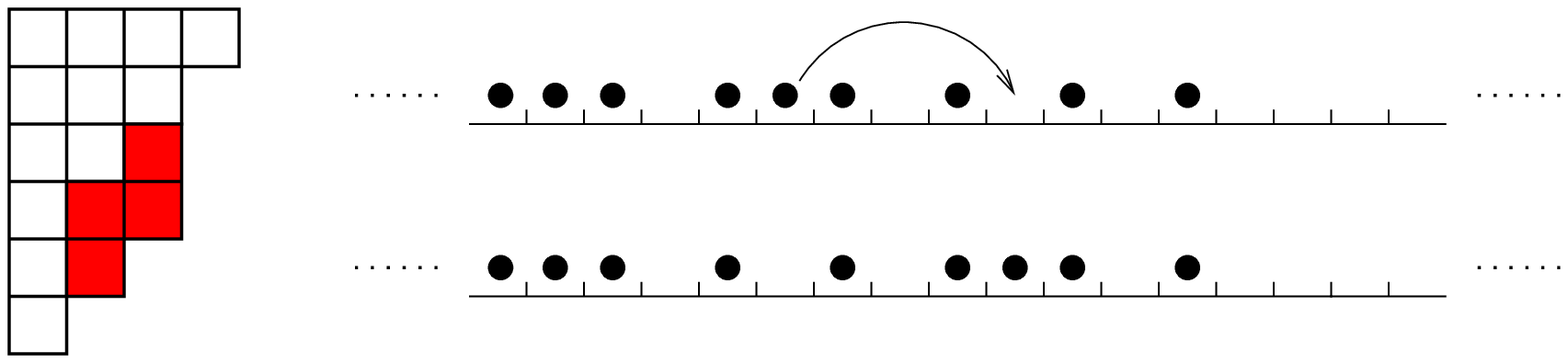}
\caption{\label{fig:coderubans} Addition of a ribbon on
$\de(\la)$.}
\end{figure}

 One may think of $\delta(\la)$ as a configuration of particles
 on the infinite discrete line: the $1$s represent particles, and
$0$s represent empty positions. So moving a particle in an empty
 position to its left (respectively right) corresponds to removing
(resp. adding) a ribbon. Figure~\ref{fig:coderubans} shows a
 ribbon $\la/\mu$ by an arrow between its tail and head, and the
codes of  $\la$ and $\mu$ are given on the right.

 One has then the following result:

\begin{lemma} \label{lem:popo}
Let $\la$ be a partition with associated sequence $\de(\la)$ ,
 and $i<j$ indices of $\de(\la)$ corresponding to a ribbon $r$
addable to (or removable from) $\la$: this just means that
 $\{\de_i(\la),\de_j(\la)\}=\{0,1\}$. Then
\begin{enumerate}
\item the size of $r$ is $|r|=j-i$.

\item the height $h$ of $r$ is the number of $1$ in $\de(\la)$
between the indices $i$ and $j$, i.e.
 $h=|\{k\in \Zm~|~i<k<j\text{ and }\de_k(\la)=1\}|$.
\end{enumerate}
\end{lemma}

\noindent\emph{Proof:} It goes simply by using the fact that the
 $1$s correspond to vertical steps on the boundary of $\la$,
and the $0$s to horizontal ones; so $j-i$ is equal to the number of cells
 occupied by $r$, and each $1$ between $i$ and $j$ corresponds
to going up from one row to another in the ribbon.\qee\medskip

Now the data of $\la,\mu,\nu$ (when $\mu,\nu\neq \la$) in a direct rule is
equivalent to the data of $\de(\la)$ and integers $i_1<j_1$, $i_2<j_2$
(corresponding to $\mu/\la$ and $\nu/\la$), where
$\de_{i_1}(\la)=1,\de_{j_1}(\la)=0$ and $\de_{i_2}(\la)=1,\de_{j_2}(\la)=0$.
Every operation of Section~\ref{subsect:oper} can be in fact easily explicated
given this representation; we shall do it for the $switchout$ operation of rule
$S$. \smallskip

  The rule $S$ applies precisely when one of the following two cases occur:
  
\begin{enumerate}
\item $i_1=i_2$, $j_1\neq j_2$ and $\de_{j_1+j_2-i_1}(\la)=1$, or
\item $j_1=j_2$, $i_1\neq i_2$ and $\de_{i_1+i_2-j_1}(\la)=0$
\end{enumerate}

\begin{figure}[!ht]
\psfrag{S}{S} \psfrag{i}{$i$} \psfrag{j}{$j_1$} \psfrag{jj}{$j_2$}
\psfrag{k}{$j_1+j_2-i$}
\includegraphics[width=\textwidth]{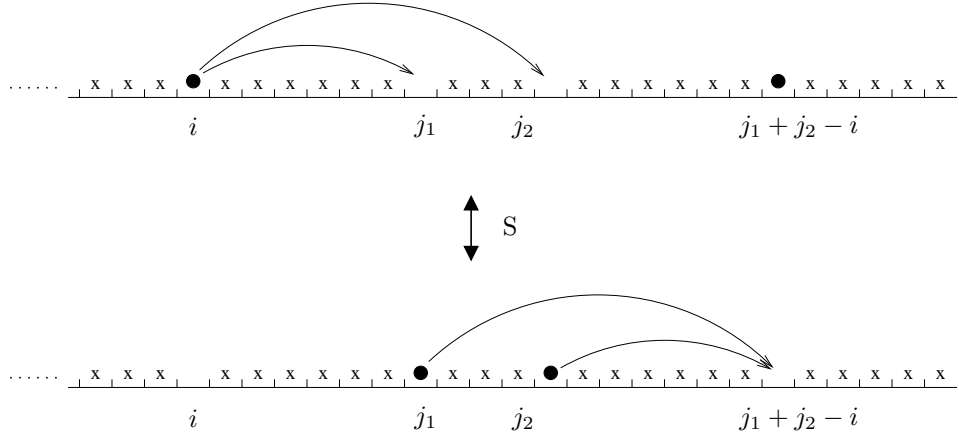}
\caption{\label{fig:srule} The rule $S$.}
\end{figure}

Let $i$ be the common value of $i_1$ and $i_2$ in the first case, and $j$ the
common value of $j_1$ and $j_2$ in the second case. Then applying rule $S$
consists simply in defining $\widehat{\la}$ as the partition whose code is
obtained from $\de(\la)$ by exchanging $0$ and $1$ at positions $i,j_1,j_2$ and
$j_1+j_2-i$ in the first case, and at positions $i_1,i_2,j$ and $i_1+i_2-j$ in
the second case. On Figure~\ref{fig:srule} the first case is illustrated: $\la$
is shown above with the ribbons, and the applications of rule $S$ gives the
partition below. (The symbols 'x' represent indifferently 1s or 0s). 

One notices that the partition $\widehat{\la}$ is element of
 $\D_{j_1-i}(\mu,\nu)$, as was $\la$: indeed, one gets the same partition
 by moving a particle following the long arrow (resp. the short arrow)
 in both partitions of Figure~\ref{fig:srule} . The rule is clearly
an involution, and in fact exchanges cases 1. and 2. defined above.

Now we will check finally that it exchanges signs, which is just a matter of
counting particles, thanks to Lemma~\ref{lem:popo}. Let $a$ (respectively $b$,
$c$) be the number of $1$ in $\la$ that are strictly between the indices $i$ and
$j_1$ (resp. $j_1$ and $j_2$, resp. $j_2$ and $j_1+j_2-i$); these numbers are
the same in $\la$ and $\widehat{\la}$. As elements of the signed set
$\D_{j_1-i}(\mu,\nu)$, $\la$ and $\widehat{\la}$ have signs $(-1)^x$ and
$(-1)^y$, with $x=2a+b$ and $y=b+2c+1$; this gives opposite signs since $x$ and
$y$ have opposite parity.
\qee

\section{Proof of Theorems~\ref{th:bij_perm} and~\ref{th:bij_invol}}
\label{sect:appa}

In this Section we will give the proofs of Theorems ~\ref{th:bij_perm} and~\ref{th:bij_invol}, the signed bijections having been defined in Section~\ref{sect:bij}. The proof is directly inspired by Fomin's constructions~\cite{FominSchen}, but some extra technicalities are needed in both proofs.
\medskip

\subsection{Proof Of Theorem~\ref{th:bij_perm}} We will show in particular that the construction $\phi$ defined algorithmically in~\ref{sub:bijhookperm} verifies indeed all properties stated in the theorem.
 We will demonstrate that this algorithm is in fact a consequence
 of Garsia and Milne's involution principle: therefore, we have
 to construct signed sets $\mathcal{A},\mathcal{B}$ and adequate
 functions. For this, a certain number of concepts have to be defined.

We call {\em border} of the grid $G_\ell$ a path from the top
left vertex to the bottom right one, with South and East steps.
 We call {\em inside} of the border $F$ the squares of $G_\ell$
 to the south west of $F$, and {\em outside} the rest of the squares.

A \emph{good labeling} of a border $F$ is the labeling of
each of its vertices by a partition, such that:

\begin{itemize}
\item the vertices at the top left and bottom right are labeled
 by the empty partition $\emptyset$.

 \item for every horizontal edge of $F$, the labels $\la$
and $\mu$ at the left and right end respectively form a ribbon $\mu/\la$.

 \item for every vertical edge of $F$, the labels $\la$
and $\mu$ at the bottom and top end respectively form a ribbon $\mu/\la$ .
\end{itemize}\medskip

Remember that we fixed a total order $\mathcal{O}$ on the squares of $G_\ell$.
Let a {\em $\mathcal{O}$-border} be a border $F$ when the squares inside $F$ are
smaller than the squares outside. Such a border defines two squares in general:
$sq_<(F)$ which is the largest square inside $F$, and $sq_>(F)$ the smallest one
outside; $sq_<(F)$ is not defined when $F$ consists of the left and bottom side
of the grid, and $sq_>(F)$ is not defined when $F$ consists of the top and right
side of the grid.

Let $F$ be a border with a good labeling $label$:
note that $label$ induces a labeling by ribbons on the edges of $F$. Suppose
 that certain squares of the grid are filled (we will also say {\em colored}) by
nonempty hooks.
 Such a coloring $col$ is \emph{compatible} with $(F,label)$ if:
\begin{itemize}
\item The squares inside $F$ are not filled.
 \item  for every horizontal edge $h$ of $F$ labeled by
 $r$, there is exactly one square filled by a hook in the column above $h$
 when $r$ is empty, and none if $r$ is non empty.
\item for every vertical edge $v$ of $F$ labeled by
 $r$, there is exactly one square filled by a hook in the row right of $v$
 when $r$ is empty, and none if $r$ is non empty.
\end{itemize}\medskip

{\bf Configurations}.
We can now introduce \textit{configurations}, which are the main objects we will consider for the rest of the proof 

\begin{definition} \label{def:conf}
 A \emph{configuration} is a 3-tuple $(F,label,col)$
 where $F$ is a $\mathcal{O}$-border, wich is well labeled by
 $label$, and $col$ is a coloring of $G_\ell$ compatible with $(F,label)$.
 \end{definition}

Let $(F,label,col)$ be a configuration on $G_\ell$. For each
 $i\in \lb 1, \ell \rb$ we note $c_i>0$ the size of the hook
in column $i$, or the size of the ribbon labeling the edge of
 $F$ appearing in column $i$: by compatibility of $label$ and
 $col$, exactly one of these two cases occur. Likewise, we
 note $c'_i$ the size of the hook in row $i$, or the size of
the ribbon labeling the edge of
 $F$ appearing row $i$. The \emph{content} of a configuration is
  then defined as the two compositions $c,c'$, where $c=(c_1,\ldots,c_\ell)$
 and $c'=(c'_1,\ldots,c'_\ell)$.

 Let us define the {\em sign} of a configuration $(F,label,col)$
as the product of all $2\ell$ ribbons labeling the edges
of $F$ (recall that the empty ribbon has sign $+1$).
For instance, the steps C and E of Figure~\ref{fig:bij}
 show two configurations : C and E have both content
$((2,1,2),(2,2,1))$, and C has negative sign, whereas
 E has positive sign.
\medskip

A hook permutation $(\s,H)$ is a positive configuration:
the border is the left and bottom side of $G_\ell$,
all vertices are labeled by $\emptyset$, and the
coloring is just the representation of $(\s,H)$
in the grid. The content of such a configuration is
 $(c(H),c(H,\s^{-1}))$.
A pair $(P,Q)$ of ribbon tableaux of length
 $\ell$ with the same shape is also a
configuration: the border is the top and right
side, the right side being labeled by $P$ and the top
 side by $Q$; and all squares of the grid are empty.
 The sign of the configuration is $\eps(P)\eps(Q)$,
and the content is $c(Q),c(Q)$. We will from now on
identify hook permutations and pairs of
ribbon tableaux to such configurations.\medskip

We now notice that in the algorithm describing our correspondence, the transformation Apply\_local\_rule 
entails a change from one configuration to another. Indeed, an
inspection of each local rule shows that the compatibility conditions in the definition
are indeed respected by this tranformation.
 \medskip

{\bf Application of the Involution principle}.
A configuration $(F,label,col)$ is {\em of type A} if
 rule $T$ has to be applied in $sq_<(F)$, or if it is a
permutation (which is when $sq_<(F)$ is not defined). It is
 {\em of type B} if rule $S$ has to be applied in
$sq_>(F)$, or if it is a pair of ribbon tableaux  (which is
when $sq_>(F)$ is not defined).

We fix two compositions $c_1$ and $c_2$ of length $\ell$ and size $n$. Let $\mathcal{A}$ (respectively  $\mathcal{B}$ ) be the configurations of type A (resp. B) and content $(c_1,c_2)$. For instance, the configurations C and E of Figure~\ref{fig:bij} are of  type $\mathcal{B}$ and $\mathcal{A}$ respectively, for the content $((2,1,2),(2,2,1))$.

Let us define an involution on $\mathcal{A}$ by applying the rule
 T in $sq_<(F)$ when it is defined, and letting the hook permutations
 unchanged.  Likewise, we have an involution on $\mathcal{B}$
   by applying the rule S in
$sq_>(F)$  when it is defined, and letting the pairs of tableaux
 unchanged. Notice that these two involutions are sign
reversing by Proposition~\ref{prop:locrul}. We define also a bijection $\mathcal{A}\rightarrow \mathcal{B}$  by applying direct rules $sq_>(F)$ on a type $\mathcal{A}$ configuration, until we reach
 a type $\mathcal{B}$ configuration. Since in this case  we only apply rules of the form Di, the sign is preserved,
thanks to Proposition~\ref{prop:locrul} again.

We now have all the functions verifying the involution principle of Garsia and Milne: this gives us a signed bijection between hook permutations and pairs of ribbon tableaux that verifies exactly the properties of Theorem \ref{th:bij_perm}, which completes the proof.\qee

\subsection{Proof of Theorem \ref{th:bij_invol}}

We will use here the definitions used in the proof of Theorem~\ref{th:bij_perm}, and adapt them to the case of involutions.
We consider a total order $\mathcal{HO}$ on $HG_\ell$, the half grid made of the squares $(i,j)$ of $G_\ell$ with $i\geq j$, and suppose that this order extends the order $\preceq$, defined by
$(i,j)\preceq (i',j')$ if $i\leq i'$ and $j\leq j'$ as seen before.
 This induces a partial order on the squares of $G_\ell$,
given by :$(i,j)$ is smaller than $(i',j')$ if and only if
$(\min(i,j),\max(i,j))$ is smaller than
$(\min(i',j'),\max(i',j'))$ for $\mathcal{HO}$. We will still note
$\mathcal{HO}$ this partial order.

\begin{figure}[!ht]
\centering
\includegraphics[width=0.5\textwidth]{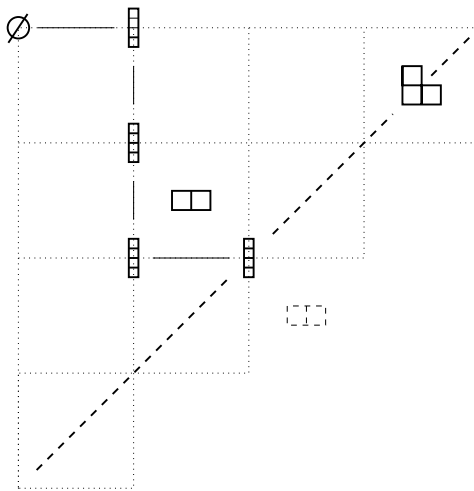}
\caption{A half configuration. 
\label{fig:growth_diagram_invol}}
\end{figure}

We consider configurations as in Definition~\ref{def:conf} ,with
the exception that the order $\mathcal{HO}$ is now used. A configuration $(F,label,col)$ on $G_\ell$ is called
 {\em symmetric} if $F,label$ and $col$ are all symmetric with
respect to the diagonal $i=j$, and we name \emph{half-configuration}
the restriction to $HG_\ell$ of a symmetric configuration.
note that the border $f$ of a half-configuration consists of a path form the top left corner of $HG_\ell$ to any vertex of the diagonal
$i=j$. Figure~\ref{fig:growth_diagram_invol} shows an example of half configuration.

For a half configuration, \emph{applying a local rule} in the
square $(i,j)\in HG_\ell$ means applying it in both $(i,j)$
and $(j,i)$ in the associated symmetric configuration, and
restrict the result to $HG_\ell$. Note that local rules in
$(i,j)$ and $(j,i)$  will give the same outputs since all
our local rules are symmetric in $\mu$  and $\nu$.
\medskip

We represent ribbon tableaux by a chain of partitions on the top side of $HG_\ell$, and hook involutions by the restriction of their matrix
representation to $HG_\ell$ with $\emptyset$ labeling the
 vertices on the left. The \emph{content} $c$ of a symmetric configuration is of the form $(c,c)$, so we define the {\em content} of the associated half configuration as $c$. The \emph{sign} of a half configuration $(f,label,col)$ is the product of all the signs labeling $f$, multiplied by the product of the signs of hooks appearing in the square $(i,i)$. Note that this gives the desired sign on $\HI$ and on ribbon tableaux, so that we are indeed in the setting of Theorem~\ref{th:bij_invol}.

Let us note $\mathcal{HA}$ and $\mathcal{HB}$ the sets of half configurations with associated symmetric configurations in $\mathcal{A}$ and $\mathcal{B}$ respectively, with the partial order $\mathcal{HO}$. We define the involutions on $\mathcal{HA}$
and $\mathcal{HB}$ in the same fashion as for $\mathcal{A}$ and $\mathcal{B}$, as well as the bijection between $\mathcal{HA}$ and $\mathcal{HB}$.
\medskip

Now we wish to apply once again the involution principle: we have here to check the sign modifications, since the definition of the sign of a configuration has been modified.

For the non diagonal squares $(i,j)$ with $i>j$, everything
works as before: the application of a rule $S$ or $T$ changes
the sign, whereas the other rules preserve it. Now there remains the application of a rule on a diagonal square $sq$. First we notice that only rules D1-D3 and I1-I3 can be applied there, since they are the only  rules with $\mu=\nu$. One needs to prove that the sign of the ribbon on the left side of $sq$, times the sign of the hook in $sq$ in the case of D2, is equal to the sign of the ribbon on the top side of $sq$.

 This is trivial for the rule D1. For the rules D2 and D3, one
 has to look at the definition of $prev$ and $first$ (Section~\ref{subsect:oper}): $first(\la,eq)$ is of the
same height as $eq$ by definition, and thus of the same sign: this implies that rule D2 will indeed have the sign preserving property. For the rule D3, $next(\mu,\mu/\la)$is a ribbon of the same height as $\mu/\la$, and so of the same sign, which implies here also that the rule preserves the sign. The involution principle can thus be applied, and this achieves the proof of Theorem~\ref{th:bij_invol}.


\begin{thebibliography}{10}

\bibitem{BesOls}
Christine Bessenrodt and J{\o}rn Olsson.
\newblock {O}n the sequence {A}085642.
\newblock {\em The {O}n-line {E}ncyclopedia of {I}nteger {S}equences}, 2004.

\bibitem{ChauveDul}
Cedric Chauve and Serge Dulucq.
\newblock A geometric version of the {R}obinson-{S}chensted correspondence for
  skew oscillating tableaux.
\newblock {\em Discrete Math.}, 246(1-3):67--81, 2002.

\bibitem{DulucqSag1}
Serge Dulucq and Bruce~E. Sagan.
\newblock The {R}obinson-{S}chensted correspondence for skew-oscillating
  tableaux.
\newblock In {\em Proceedings of the 4th conference on Formal power series and
  algebraic combinatorics}, pages 129--142, Amsterdam, The Netherlands, The
  Netherlands, 1995. Elsevier Science Publishers B. V.

\bibitem{DulucqSag2}
Serge Dulucq and Bruce~E. Sagan.
\newblock La correspondance de {R}obinson-{S}chensted pour les tableaux oscillants gauches.
\newblock {\em Discrete Mathematics}, 139(1-3):129-142, 1995.

\bibitem{FominGen}
S.~V. Fomin.
\newblock The generalized {R}obinson-{S}chensted-{K}nuth correspondence.
\newblock {\em Zap. Nauchn. Sem. Leningrad. Otdel. Mat. Inst. Steklov. (LOMI)},
  155(Differentsialnaya Geometriya, Gruppy Li i Mekh. VIII):156--175, 195,
  1986.

\bibitem{FominDual}
Sergey Fomin.
\newblock Duality of graded graphs.
\newblock {\em J. Algebraic Combin.}, 3(4):357--404, 1994.

\bibitem{FominSchen}
Sergey Fomin.
\newblock Schensted algorithms for dual graded graphs.
\newblock {\em J. Algebraic Combin.}, 4(1):5--45, 1995.

\bibitem{FominSchur}
Sergey Fomin.
\newblock Schur operators and {K}nuth correspondences.
\newblock {\em J. Combin. Theory Ser. A}, 72(2):277--292, 1995.

\bibitem{FominStan}
Sergey Fomin and Dennis Stanton.
\newblock Rim hook lattices.
\newblock {\em Algebra i Analiz}, 9(5):140--150, 1997.

\bibitem{FultonBook}
William Fulton.
\newblock {\em Young tableaux}, volume~35 of {\em London Mathematical Society
  Student Texts}.
\newblock Cambridge University Press, Cambridge, 1997.
\newblock With applications to representation theory and geometry.

\bibitem{GarMilshort}
Adriano Garsia and Stephen Milne.
\newblock Method for constructing bijections for classical partition
  identities.
\newblock {\em Proc. Nat. Acad. Sci. U.S.A.}, 78(4, part 1):2026--2028, 1981.

\bibitem{GarMil}
Adriano Garsia and Stephen Milne.
\newblock A {R}ogers-{R}amanujan bijection.
\newblock {\em J. Combin. Theory Ser. A}, 31(3):289--339, 1981.

\bibitem{Isaacs}
I.~Martin Isaacs.
\newblock {\em Character theory of finite groups}.
\newblock Dover Publications Inc., New York, 1994.
\newblock Corrected reprint of the 1976 original [Academic Press, New York;
  MR0460423 (57 \#417)].

\bibitem{Kerber}
Adalbert Kerber.
\newblock {\em Applied finite group actions}, volume~19 of {\em Algorithms and
  Combinatorics}.
\newblock Springer-Verlag, Berlin, second edition, 1999.

\bibitem{Macdonald}
I.~G. Macdonald.
\newblock {\em Symmetric functions and orthogonal polynomials}, volume~12 of
  {\em University Lecture Series}.
\newblock American Mathematical Society, Providence, RI, 1998.
\newblock Dean Jacqueline B. Lewis Memorial Lectures presented at Rutgers
  University, New Brunswick, NJ.

\bibitem{Mur}
F.~D. Murnaghan.
\newblock On the {R}epresentations of the {S}ymmetric {G}roup.
\newblock {\em Amer. J. Math.}, 59(3):437--488, 1937.

\bibitem{Nak}
T.~Nakayama.
\newblock On some modular properties of irreducible representations of a
  symmetric group. {I}.
\newblock {\em Jap. J. Math.}, 18:89--108, 1941.

\bibitem{Roby}
Tom Roby.
\newblock The connection between the {R}obinson-{S}chensted correspondence for
  skew oscillating tableaux and graded graphs.
\newblock {\em Discrete Math.}, 139(1-3):481--485, 1995.
\newblock Formal power series and algebraic combinatorics (Montreal, PQ, 1992).

\bibitem{SaganBook}
Bruce~E. Sagan.
\newblock {\em The symmetric group}, volume 203 of {\em Graduate Texts in
  Mathematics}.
\newblock Springer-Verlag, New York, second edition, 2001.
\newblock Representations, combinatorial algorithms, and symmetric functions.

\bibitem{SaganStanley}
Bruce~E. Sagan and Richard P. Stanley.
\newblock {R}obinson-{S}chensted algorithms for {S}kew tableaux.
\newblock {\em J. Comb. Theory Ser. A}, 55(2):161--193, 1990.

\bibitem{Sagan}
Bruce~E. Sagan.
\newblock Shifted tableaux, {S}chur {Q}-functions, and a conjecture of {R}. {S}tanley.
\newblock {\em J. Comb. Theory Ser. A}, 45(1):62--103, 1987.

\bibitem{Schen}
C.~Schensted.
\newblock Longest increasing and decreasing subsequences.
\newblock {\em Canad. J. Math.}, 13:179--191, 1961.

\bibitem{ShiWhi}
Mark Shimozono and Dennis~E. White.
\newblock Color-to-spin ribbon {S}chensted algorithms.
\newblock {\em Discrete Math.}, 246(1-3):295--316, 2002.
\newblock Formal power series and algebraic combinatorics (Barcelona, 1999).

\bibitem{Online}
Neil Sloane.
\newblock On-line encyclopedia of integer sequences.
\newblock Accessible from N. Sloane's homepage.

\bibitem{StanleyDiff}
Richard~P. Stanley.
\newblock Differential posets.
\newblock {\em J. Amer. Math. Soc.}, 1(4):919--961, 1988.

\bibitem{StanleyVarDiff}
Richard~P. Stanley.
\newblock Variations on differential posets.
\newblock In {\em Invariant theory and tableaux (Minneapolis, MN, 1988)},
  volume~19 of {\em IMA Vol. Math. Appl.}, pages 145--165. Springer, New York,
  1990.

\bibitem{StanleyEnum2}
Richard~P. Stanley.
\newblock {\em Enumerative combinatorics. {V}ol. 2}, volume~62 of {\em
  Cambridge Studies in Advanced Mathematics}.
\newblock Cambridge University Press, Cambridge, 1999.
\newblock With a foreword by Gian-Carlo Rota and appendix 1 by Sergey Fomin.

\bibitem{StanWhi}
Dennis~W. Stanton and Dennis~E. White.
\newblock A {S}chensted algorithm for rim hook tableaux.
\newblock {\em J. Combin. Theory Ser. A}, 40(2):211--247, 1985.

\bibitem{Sundaram}
Sheila Sundaram.
\newblock The {C}auchy identity for {$\rm Sp(2n)$}
\newblock {\em J. Combin. Theory Ser. A}, 53(2):209--238, 1990.

\bibitem{VanLeeuEdge}
Marc A.~A. van Leeuwen.
\newblock Edge sequences, ribbon tableaux, and an action of affine
  permutations.
\newblock {\em European J. Combin.}, 20(2):179--195, 1999.

\bibitem{WhiteBij}
Dennis~E. White.
\newblock A bijection proving orthogonality of the characters of {$S\sb{n}$}.
\newblock {\em Adv. in Math.}, 50(2):160--186, 1983.

\end{thebibliography}
\end{document}